\documentclass[11pt]{article}
\usepackage{amsmath, graphicx, amsfonts,amssymb, calrsfs}
\usepackage{amsfonts,mathrsfs, color, amsthm}
 
\usepackage{enumitem}
 
\def\cK{\mathcal{K}}
\def\sphere{S^{n-1}}

\def\Rn{{\mathbb R^n}}
 \def\R{\mathbb{R}}
\def\cH{\mathcal{H}}

\newtheorem{theorem}{Theorem}[section]
\newtheorem{lemma}{Lemma}[section]
\newtheorem{remark}{Remark}[section]
\newtheorem{proposition}{Proposition}[section]
\newtheorem{corollary}{Corollary}[section]

\newtheorem{definition}{Definition}[section]

\def\cC{\mathcal{C}}
\def\bpf{\begin{proof}}
\def\epf{\end{proof}}
\def\be{\begin{equation}}
\def\ee{\end{equation}}
\def\bea{\begin{eqnarray}}
\def\eea{\end{eqnarray}}
\def\bt{\begin{theorem}}
\def\et{\end{theorem}}
\def\bl{\begin{lemma}}
\def\el{\end{lemma}}
\def\br{\begin{remark}}
\def\er{\end{remark}}
\def\bc{\begin{corollary}}
\def\ec{\end{corollary}}
\def\bd{\begin{definition}}
\def\ed{\end{definition}}
\def\bp{\begin{proposition}}
\def\ep{\end{proposition}}

\begin{document}
\title{The $p$-capacitary Orlicz-Hadamard variational formula and Orlicz-Minkowski problems  }

\author{Han Hong, Deping Ye and Ning Zhang}
\date{}
\maketitle
\begin{abstract}  

In this paper, combining the $p$-capacity for $p\in (1, n)$ with the Orlicz addition of convex domains, we develop the $p$-capacitary Orlicz-Brunn-Minkowski theory. In particular, the Orlicz $L_{\phi}$ mixed $p$-capacity of two convex domains is introduced and its geometric interpretation is obtained by the $p$-capacitary Orlicz-Hadamard variational formula. The $p$-capacitary Orlicz-Brunn-Minkowski and Orlicz-Minkowski inequalities are established, and the equivalence of these two inequalities are discussed as well. The $p$-capacitary Orlicz-Minkowski problem is proposed and  solved under some mild conditions on the involving functions and measures.  In particular, we provide the solutions for the normalized $p$-capacitary $L_q$ Minkowski problems with $q>1$ for both discrete and general measures. 

 \end{abstract}\section{Introduction}
 
 The classical Minkowski problem aims to find the necessary and/or sufficient conditions on a given finite Borel measure $\mu$ defined on the unit sphere $\sphere\subset \Rn$ such that $\mu$ is the surface area measure of a convex body (i.e., a convex and compact subset of $\Rn$ with nonempty interior). Its $L_q$ extension, namely the $L_q$ Minkowski problem \cite{Lut1993},  has been a central object of interest in convex geometric analysis for decades and has received extensive considerations (see e.g., \cite{CW, CWXJ, HuMaShen,  HugLYZ, LYZ2004, Umanskiy,  ZhuGX2015a,  ZhuGX2015b, ZhuGX2017}). Both the classical and $L_q$ Minkowski problems are related to  function $\varphi=t^q$ for $0\neq q\in \R$. There are versions of Minkowski problems related to other functions, for instance, the $L_0$ Minkowski or logarithmic Minkowski problems \cite{BorHegZhu, bor2013-2, Stancu2002, Stancu2003, Stancu2008,  ZhuGX2014} and the Orlicz-Minkowski problem \cite{HLYZ, huang}.

Replacing the surface area measure in the classical Minkowski problem by the $p$-capacitary measure for  $p\in (1, n)$, the following $p$-capacitary $L_1$ Minkowski problem can be asked and is of central importance in the development of the $p$-capacitary Brunn-Minkowski theory: {\em under what conditions on a given finite Borel measure $\mu$ defined on $\sphere$, one can find a convex domain (i.e., the interior of a convex body) whose $p$-capacitary measure is equal to $\mu$?}  When $p=2$, this has been solved in the seminal papers by Jerison \cite{Jerison, Jerison-1996}. A solution of this problem for  $p\in (1, 2)$ was given  by Colesanti, Nystr\"om, Salani, Xiao, Yang and Zhang in their remarkable paper \cite{CNSXYZ}. {The  normalized (nonlinear) $p$-capacitary  $L_1$ Minkowski problem} for all $p\in (1, n)$ was recently solved by Akman, Gong, Hineman, Lewis and Vogel in their groundbreaking paper \cite{AGHLV}, where the necessary and sufficient conditions for the finite Borel measure $\mu$ being the $p$-capacitary measure of a convex domain were provided.

In view of the classical Minkowski problem and its various extensions, it is important to investigate the $p$-capacitary $L_q$ Minkowski and Orlicz-Minkowski problems. More precisely, we propose the following question: {\em under what conditions on a given function $\phi$ and a given finite Borel measure $\mu$ defined on $\sphere$, one can find a convex domain $\Omega$ such that the origin $o$ is in its closure and 
$$
   \frac{\mu}{\phi(h_{\Omega})}=\tau \cdot \mu_p(\Omega,\cdot),
$$  
where $\tau>0$ is a constant?} Hereafter, $\mu_p(\Omega, \cdot)$ defined on $\sphere$ denotes the $p$-capacitary measure of $\Omega$, and  $h_{\Omega}$ denotes the support function of $\Omega$ (see Section \ref{section 2} for notations). When $\phi=t^{q-1}$ for $q\in \R$, one gets the following  normalized $p$-capacitary  $L_q$ Minkowski problem: {\em under what conditions on a given finite Borel measure $\mu$ defined on $\sphere$, one can find a convex domain $\Omega$ such that the origin $o$ is in its closure and 
$$
   \mu \cdot h^{q-1}_{\Omega}\cdot C_p(\Omega)=c(n, p, q, \mu) \cdot  \mu_p(\Omega,\cdot), 
$$ 
where $c(n, p, q, \mu)>0$ is a constant and $C_p(\Omega)$ is the $p$-capacity of $\Omega$?}  In Subsection \ref{section:polytope-minkowski-solution}, we provide a solution for the above $p$-capacitary Minkowski problems for discrete measures under some very limited assumptions on $\mu$: {\em the support of $\mu$ is not contained in any closed hemisphere}.  A solution of the above $p$-capacitary Minkowski problems for general measures is provided in Subsection \ref{M-capacity-general--111}.   

The $p$-capacitary measure can be derived from an integral related to the $p$-equilibrium potential of $\Omega$. Note that the $p$-equilibrium potential of $\Omega$ is the solution of a $p$-Laplace equation with certain boundary conditions (see Subsection \ref{Section: p-capacity-1} for details). For a convex domain $\Omega\subset \Rn$, the Poincar\'e $p$-capacity formula  \cite{CNSXYZ} asserts that the $p$-capacity of $\Omega$ has the following form:  
\begin{equation} \label{calculation-p-capacity} 
   C_p(\Omega)=\frac{p-1}{n-p}\int_{S^{n-1}} h_\Omega(u)\,d\mu_p(\Omega, u).
\end{equation}  
Although the definition of the $p$-capacity involves rather complicate partial differential equations, formula (\ref{calculation-p-capacity}) suggests that the $p$-capacity has high resemblance with the volume. For a convex domain $\Omega\subset \Rn$, its volume can be calculated by  
$$
   |\Omega|=\frac{1}{n} \int _{\sphere}  h_\Omega(u)\,dS(\Omega, u),
$$ 
with $S(\Omega, \cdot)$ the surface area measure of $\Omega$ defined on $\sphere$.   For instance, the $p$-capacitary Brunn-Minkowski inequality asserts that for all convex domains $\Omega$ and $\Omega_1$, one has 
\begin{equation}
   \label{class-B-M-capacity-20160615} C_p(\Omega + \Omega_1)^{\frac{1}{n-p}}\geq C_p(\Omega)^{\frac{1}{n-p}}+C_p(\Omega_1)^{\frac{1}{n-p}},
\end{equation}  
with equality if and only if $\Omega$ and  $\Omega_1$ are homothetic (see \cite{Borell, Caffarelli1996, CS2003}).   Hereafter 
$$
   \Omega + \Omega_1=\{x+y: x\in \Omega, y\in \Omega_1\}
$$
denotes the Minkowski sum of $\Omega$ and $\Omega_1$. Inequality (\ref{class-B-M-capacity-20160615}) is  similar to the classical Brunn-Minkowski inequality regarding the volume: \begin{equation*}  
   |\Omega + \Omega_1|^{\frac{1}{n}}\geq |\Omega|^{\frac{1}{n}}+|\Omega_1|^{\frac{1}{n}},
\end{equation*}  
with equality if and only if $\Omega$ and  $\Omega_1$ are homothetic (see e.g., \cite{GardnerB-M,SchneiderBook}). Moreover, the $p$-capacitary Minkowski inequality (\ref{Minkowski-classical-1}) shares the formula similar to its volume counterpart (see e.g., \cite{CNSXYZ, GardnerB-M, SchneiderBook}).

Sections \ref{Section:OM-Inequality} and \ref{Section:OBM} in this paper reveal  another surprising similarity between the $p$-capacity and the volume regarding the Orlicz additions. We develop the $p$-capacitary Orlicz-Brunn-Minkowski theory based on the combination of the Orlicz additions and the $p$-capacity.   The Orlicz additions were proposed by Gardner, Hug and Weil  in  \cite{GHW2013} and independently by Xi, Jin  and Leng in \cite{XJL}, in order to provide the foundation of  the newly initiated Orlicz-Brunn-Minkowski theory for convex bodies (with respect to volume) starting from the works \cite{LYZ2010a, LYZ2010b} of Lutwak, Yang and Zhang. The Orlicz theory is in great demand (see e.g., \cite{Ye2016} for some motivations) and is rapidly developing (see e.g., \cite{bor2013, BLYZ, HaberlParapatits2014, Ludwig2009, XL2016,  Ye2012, Ye2013, Ye2014, ZouXiong2014}).  In particular, we establish the $p$-capacitary Orlicz-Brunn-Minkowski inequality (see Theorem \ref{C-O-B-M-inequality-1}) and Orlicz-Minkowski inequality (see Theorem \ref{capacitary-Minkowski-2-2-1-1}). The $p$-capacitary Orlicz-Minkowski inequality provides a tight lower bound for $C_{p, \phi}(\Omega, \Omega_1)$, the Orlicz $L_{\phi}$ mixed $p$-capacity of $\Omega, \Omega_1\in \cC_0$ (the collection of all convex domains containing the origin),  in terms of $C_p(\Omega)$ and $C_p(\Omega_1)$.  In Theorem \ref{c-interpretation-111}, we prove the $p$-capacitary Orlicz-Hadamard variational formula based on a linear Orlicz addition of $\Omega, \Omega_1\in \cC_0$. This $p$-capacitary Orlicz-Hadamard variational formula gives a geometric interpretation of $C_{p, \phi}(\Omega, \Omega_1)$.  Section \ref{section 2} is for the necessary background and notation. More details could be found in \cite{CNSXYZ, EvansGariepy,  SchneiderBook}.  It is worth to mention that results in this paper are for the $p$-capacity related to the $p$-Laplace equations; however similar results for the nonlinear $\mathcal{A}$-capacity associated with a nonlinear elliptic partial differential equation \cite{AGHLV} could be obtained as well.

\section{Background and Notations}\label{section 2}

Throughout this paper, $n\geq 2$ is a natural number. A subset $K\subset\Rn$ is convex if $\lambda x+(1-\lambda)y\in K$ for any $x, y\in K$ and $\lambda\in [0, 1]$. 
A convex set $K\subset \Rn$ is a convex body if $K$ is also compact with nonempty interior.  Denote by $\cK_0$ the set of convex bodies in $\Rn$ with the origin in their interiors.  The usual Euclidean norm is written by $\|\cdot\|$ and the origin of $\Rn$ is denoted by $o$.  Let $\{e_1,\cdots, e_n\}$ be the standard orthonormal basis of $\Rn$. 
Define $\lambda K=\{\lambda x: x\in K\}$ for $\lambda\in \R$ and $K\subset \Rn$. For a convex body $K\in \cK_0$,  $|K|$  refers to the volume of $K$ and $\cH^{n-1}$ refers to the $n-1$ dimensional Hausdorff measure of $\partial K$, the boundary of $K$. For a set $E\subset \R^n$, define $conv(E)$, the convex hull of $E$, to be the smallest convex set containing $E$.  Let $\theta^{\perp}=\{x\in \Rn: \langle x, \theta\rangle=0\}$ for $\theta\in \sphere$.

The support function of a convex compact set $K$ containing the origin is the function $h_K: S^{n-1}\rightarrow [0, \infty)$ defined by 
$$
   h_K(u)=\max_{y\in K}\langle y, u\rangle,
$$ 
where $\langle \cdot, \cdot\rangle$ is the usual inner product on $\Rn$.  Hereafter, $S^{n-1}$ is the unit sphere of $\Rn$ which consists of all unit vectors in $\Rn$. Note that the support function $h_K$  can be extended to $\Rn\setminus \{o\}$ by 
$$
   h_K(x)=h_K(ru)=rh_K(u)
$$ 
for $x=ru$ with $u\in S^{n-1}$ and $r\geq 0$. Clearly $h_K: \sphere \rightarrow \R$ is sublinear.  Any convex body $K\in \cK_0$ is uniquely characterized by its support function.  Two convex bodies $K, L\in \cK_0$ are said to be dilates of each other if $h_K=c\cdot h_L$ for some constant $c>0$; $K$ and $L$ are said to be homothetic to each other if $K-a$  dilates  $L$ for some $a\in \Rn$.  On $\cK_0$, the Hausdorff metric $d_H(\cdot, \cdot)$ is a natural way to measure the distance of two convex bodies $K, L\in \cK_0$: 
$$
   d_H(K, L)=\max_{u\in \sphere} |h_K(u)-h_L(u)|= \|h_K-h_L\|_{\infty}. 
$$
The Blaschke selection theorem (see e.g., \cite{SchneiderBook}) states that  {\em  every bounded sequence of convex bodies has a subsequence
that converges to a (possibly degenerated) convex compact set.}

 Note that $K\in \cK_0$ can be formulated by  the intersection of hyperspaces as follows:
$$
   K=\mathop{\bigcap}\limits_{u\in S^{n-1}} \Big\{x\in \Rn: \langle x, u\rangle \leq h_K(u)\Big\}.
$$  
By  $C^{+}(S^{n-1})$ we mean the set of all continuous and positive functions defined on $S^{n-1}$.  The metric $d(\cdot, \cdot)$ on $C^{+}(S^{n-1})$ is assumed to be the one induced by the maximal norm: for all $f, g\in C^{+}(S^{n-1})$,  
$$
   d(f, g)=\|f-g\|_{\infty} =\max_{u\in \sphere}  |f(u)-g(u)|.
$$  
Associated to each $f\in  C^{+}(S^{n-1})$, one can define a convex body $K_f\in \cK_0$ (see formula (7.97) in \cite{SchneiderBook}) by  
$$
   K_f=\mathop{\bigcap}\limits_{u\in S^{n-1}} \Big\{x\in \Rn: \langle x,  u\rangle \leq f(u) \Big\}.
$$ 
The convex body $K_f$ is called  the  Aleksandrov body  associated to  $f \in C^{+}(S^{n-1}).$ The  Aleksandrov body provides a powerful tool in convex geometry and plays crucial roles in this paper. Here we list some important properties for  the Aleksandrov body which will be used in later context. These properties and the proofs can be found in section 7.5 in \cite{SchneiderBook}.  First of all,  if $f\in C^+(\sphere)$ is the support function of a convex body $K\in \cK_0$, then $K=K_f$. Secondly,  for $f\in C^+(\sphere)$,  $h_{K_f}(u)\leq f(u)$ for all $u\in \sphere$,  and $h_{K_f}(u)=f(u)$  almost everywhere with   respect   to  $S(K_f, \cdot)$, the surface area measure of $K_f$ defined on $\sphere$ (see the proof of Lemma 7.5.1 in \cite{SchneiderBook}). Recall that  $S(K, \cdot)$ has the following geometric interpretation  (see e.g., \cite[page 111]{SchneiderBook}):  for any Borel set $\Sigma \subset S^{n-1}$, 
\begin{equation}\label{Def:surface area----123} 
   S(K, \Sigma )=\cH^{n-1}\{x\in \partial K:  \mathrm{g}(x)\in \Sigma\},
\end{equation}  
 where $\mathrm{g}: \partial K \rightarrow S^{n-1}$ is the (single-valued) Gauss map of $K$, that is, $\mathrm{g}(x)\in S^{n-1}$ is the unit outer normal vector of $\partial K$ at almost everywhere $x\in \partial K$ with respect to the $(n-1)$-dimensional Hausdorff measure of $\partial K$. Furthermore, the convergence of $\{K_{f_m}\}_{m\geq 1}$ in the Hausdorff metric is guaranteed by the convergence of $\{f_m\}_{m\geq 1}.$ This is the Aleksandrov's convergence lemma \cite{Aleks1938} (see also  \cite[Lemma 7.5.2]{SchneiderBook}):  if the sequence $f_1, f_2, \cdots  \in C^{+}(S^{n-1})$ converges to $f\in C^{+}(S^{n-1})$ in the metric $d(\cdot, \cdot)$, then $ K_{f_1}, K_{f_2}, \cdots \in \cK_0$ converges to  $K_f\in \cK_0$ with respect to the Hausdorff metric. For more background on convex geometry, please refer to \cite{SchneiderBook}.

 \subsection{Orlicz addition  and the Orlicz-Brunn-Minkowski theory of convex bodies}   \label{section:OBM}
Let $m\geq 1$ be an integer number. Denote by $\Phi_m$ the set of convex functions $\varphi: [0,\infty)^m\to [0,\infty)$ that are increasing in each variable, and satisfy $\varphi (o)=0$ and $\varphi(e_j)=1$ for $j=1,\dots,m$.  The Orlicz $L_\varphi$ sum of $K_1, \cdots, K_m\in \cK_0$ \cite{GHW2014} is the convex body $+_{\varphi}(K_1,\dots,K_m)$ whose support function $h_{+_{\varphi}(K_1,\dots,K_m)}$ is defined by the unique positive solution of the following equation:  
\begin{equation*} 
   \varphi\left( \frac{h_{K_1}(u)}{\lambda},\dots, \frac{h_{K_m}(u)}{\lambda}\right) = 1, \ \ \ \ \mathrm{for}\ \ \  u\in \sphere. 
\end{equation*}   
That is, for each fixed $u\in S^{n-1}$,  
\begin{equation*} 
   \varphi\bigg( \frac{h_{K_1}(u)}{h_{+_{\varphi}(K_1,\dots, K_m)}(u)},\dots, \frac{h_{K_m}(u)}{h_{+_{\varphi}(K_1,\dots, K_m)}(u)}\bigg) = 1.
\end{equation*}  
The fact that $\varphi\in \Phi_m$ is increasing in each variable implies that, for $j=1, \cdots, m$,  
\begin{equation}\label{compare-Orliczdef} 
   K_j\subset +_{\varphi}(K_1,\dots,K_m). 
\end{equation} 
It is easily checked that if $K_i$  for all $1<i\leq m$ are dilates of $K_1$,  then $ +_{\varphi}(K_1,\dots,K_m)$ is dilate of $K_1$ as well.  The related Orlicz-Brunn-Minkowski inequality has the following form \cite{GHW2014}: 
\begin{equation}\label{Orldef-001}
   \varphi\bigg(\frac{|K_1|^{1/n}}{|+_{\varphi}(K_1,\dots, K_m)|^{1/n}},\dots, \frac{|K_m|^{1/n}}{|+_{\varphi}(K_1,\dots, K_m)|^{1/n}}\bigg) \leq 1.
\end{equation}  
The classical  Brunn-Minkowski and the $L_q$ Brunn-Minkowski inequalities are associated to $\varphi(x_1, \cdots, x_m)=\sum_{i=1}^mx_i\in \Phi_m$ and $\varphi(x_1, \cdots, x_m)=\sum_{i=1}^m x_i^q\in \Phi_m$ with $q>1$, respectively.  In these cases, the $L_q$ sum of $K_1,\cdots, K_m$ for $q\geq 1$ is the convex body $K_1+_q  \cdots +_qK_m$ whose support function is formulated by 
$$
   h_{K_1+_q  \cdots +_qK_m}^q=h_{K_1}^q+\cdots+h_{K_m}^q.
$$ 
When $q=1$, we often write $K_1+\cdots +K_m$ instead of $K_1+_1\cdots +_1K_m$. 

Consider the convex body $K+_{\varphi, \varepsilon}L\in \cK_0$ whose support function is given by, for $u\in \sphere$,  \begin{equation}\label{definition-linear-addition} 
   1= \varphi_1\left( \frac{h_K(u)}{h_{K+_{\varphi, \varepsilon}L}(u)}\right)+\varepsilon \varphi_2\left( \frac{h_L(u)}{h_{K+_{\varphi, \varepsilon}L}(u)}\right),
\end{equation}  
where $\varepsilon>0$, $K, L\in \cK_0$, and $\varphi_1, \varphi_2\in \Phi_1$.  If $(\varphi_1)'_l(1)$, the left derivative of $\varphi_1$ at $t=1$, exists and is positive, then  the $L_{\varphi_2}$ mixed volume of $K, L\in \cK_0$ can be defined by \cite{GHW2014, XJL, ZHY2016}
\begin{equation}\label{geometricinterpretation-mixedvolume}
   V_{\varphi_2}(K, L)=\frac{(\varphi_1)'_l(1)}{n}
   \cdot  \frac{\,d }{\,d\varepsilon}|K+_{\varphi, \varepsilon}L|
   \bigg|_{\varepsilon=0^+}= \frac{1}{n} \int_{S^{n-1}}
   \varphi_2\bigg(\frac{h_L(u)}{h_K(u)}\bigg)h_K(u)dS(K, u).
\end{equation} 
Together with the Orlicz-Brunn-Minkowski inequality (\ref{Orldef-001}), one gets the following fundamental Orlicz-Minkowski inequality: if $\varphi\in \Phi_1$, then for all $K, L\in \cK_0$, $$V_{\varphi}(K, L)\geq |K|\cdot \varphi\bigg(\bigg(\frac{|L| }{|K|}\bigg)^{1/n}\bigg),$$ with equality, if in addition $\varphi$ is strictly convex,  if and only if $K$ and $L$ are dilates of each other. The classical Minkowski and the $L_q$ Minkowski inequalities are associated with $\varphi=t$ and $\varphi=t^q$ for $q>1$ respectively.

Formula (\ref{geometricinterpretation-mixedvolume}) was proved in \cite{GHW2014, XJL} with assumptions  $\varphi_1, \varphi_2\in \Phi_1$ (i.e.,  convex and increasing functions);  however, it can be extended to more general increasing or decreasing functions  \cite{ZHY2016}. To this end, we work on the following classes of nonnegative continuous functions:  
\begin{eqnarray*} 
   \mathcal{I} \!\!\!&=&\!\!\! \{\phi: [0, \infty)\rightarrow [0, \infty) \ \mbox{such that} \ \phi \ \mbox{is strictly increasing with}\ \phi(1)=1, \phi(0)=0, \phi(\infty)=\infty\}, \\ \mathcal{D} \!\!\! &=&\! \!\! \{\phi:  (0, \infty)\rightarrow (0, \infty)\ \mbox{such that} \ \phi \ \mbox{is strictly decreasing with} \ \phi(1)=1, \phi(0)=\infty, \phi(\infty)=0\},
\end{eqnarray*}     
where for simplicity we let $\phi(0)=\lim_{t\rightarrow 0^+} \phi(t)$ and $\phi(\infty)=\lim_{t\rightarrow \infty} \phi(t)$. Note that results may still hold if the normalization on $\phi(0), \phi(1)$  and $\phi(\infty)$ are replaced by other quantities. The linear Orlicz addition of $h_K$ and $h_L$  in formula (\ref{definition-linear-addition})  can be defined in the same way for  either $\varphi_1, \varphi_2\in \mathcal{I}$ or $\varphi_1, \varphi_2\in \mathcal{D}$. Namely,  for  either $\varphi_1, \varphi_2\in \mathcal{I}$ or $\varphi_1, \varphi_2\in \mathcal{D}$, and for  $\varepsilon>0$, define $f_{\varepsilon}: \sphere\rightarrow (0, \infty)$  the linear Orlicz addition of $h_K$ and $h_L$ by, for $u\in \sphere$, 
\begin{equation}\label{definition-linear-orlicz-extension-1}
   \varphi_1\left(\frac{h_K(u)}{f_{\varepsilon}(u)}\right)+\varepsilon\varphi_2
   \left(\frac{h_L(u)}{f_{\varepsilon} (u)}\right)=1.
\end{equation}
See  \cite{HouYe2016} for more details. In general, $f_{\varepsilon}$  may not be the support function of a convex body; however $f_{\varepsilon}$ is the support function of $K+_{\varphi, \varepsilon}L$  when $\varphi_1,\varphi_2\in \Phi_1$.   It is easily checked that
$f_\varepsilon\in C^+(\sphere)$  for all $\varepsilon>0$. Moreover, $h_K\leq f_\varepsilon$ if
$\varphi_1,\varphi_2\in \mathcal{I}$ and $h_K\geq f_\varepsilon$ if
$\varphi_1,\varphi_2\in \mathcal{D}$. Denote by $K_{\varepsilon}$ the Aleksandrov body associated to $f_{\varepsilon}$.  The following result \cite{ZHY2016} extends formula (\ref{geometricinterpretation-mixedvolume}) to not necessarily convex functions $\varphi_1$ and $\varphi_2$: if $K, L\in\mathcal{K}_0$ and $\varphi_1,\varphi_2\in \mathcal{I}$ are such that $(\varphi_1)'_l(1)$ exists and is positive, then 
\begin{equation}\label{convergence of f-e---1}
   V_{\varphi_2}(K, L)=\frac{(\varphi_1)'_l(1)}{n}
   \cdot  \frac{\,d }{\,d\varepsilon}|K_{\varepsilon}|
   \bigg|_{\varepsilon=0^+}= \frac{1}{n} \int_{S^{n-1}}
   \varphi_2\bigg(\frac{h_L(u)}{h_K(u)}\bigg)h_K(u)dS(K, u),
\end{equation} 
while if $\varphi_1,\varphi_2\in \mathcal{D}$ such that $(\varphi_1)'_r(1)$, the right derivative of $\varphi_1$ at $t=1$,  exists and is nonzero,  then  (\ref{convergence of f-e---1}) holds with $(\varphi_1)'_l(1)$ replaced by $(\varphi_1)'_r(1)$.

\subsection{The $p$-capacity} \label{Section: p-capacity-1}
Throughout this paper, the standard notation  $C_c^{\infty} (\Rn)$ denotes the set of all infinitely differentiable functions with compact support in $\Rn$ and   $\nabla f$ denotes the gradient of $f$.   Let $n\geq 2$ be an integer and  $p\in (1, n)$. The $p$-capacity of a compact subset  $E\subset \Rn$, denoted by $C_p(E)$, is defined by 
$$
   C_p(E)=\inf\left\{\int_{\Rn} \|\nabla f\|^p\,dx: f\in C_c^{\infty}(\Rn) \ \ \mathrm{such\ that}\ \ f\geq 1\ \ \mathrm{on}\ \ E \right\}. 
$$  
If $O\subset \Rn$ is an open set, then the $p$-capacity of $O$ is defined by 
$$
   C_p(O)=\sup \big\{C_p(E):  \ \ E\subset O \ \ \mbox{and $E$ is a compact set in $\Rn$} \big\}.
$$ 
For general bounded measurable subset $F\subset \Rn$, the $p$-capacity of $F$ is then defined by 
$$
   C_p(F) =\inf \big\{C_p(O):  \ \ F\subset O \ \ \mbox{and $O$ is an open set in $\Rn$} \big\}.
$$ 

The $p$-capacity is monotone, that is, if $A\subset B$ are two measurable subsets of $\Rn$, then $C_p(A)\leq C_p(B)$. It is translation invariant: $C_p(F+x_0)=C_p(F)$ for all $x_0\in \Rn$ and measurable subset $F\subset\Rn$.  Its homogeneous degree is $n-p$, i.e.,  for all $\lambda>0$,  
\begin{equation} \label{homogeneous-111111111111} 
    C_p(\lambda A)=\lambda^{n-p}C_p(A).
\end{equation}  
For $K\in \cK_0$,  let $int(K)$ denote its interior. It follows from the monotonicity of the $p$-capacity that $C_p(int(K))\leq C_p(K)$.  On the other hand, for all $\varepsilon>0$, one sees that $$K\subset (1+\varepsilon)\cdot int(K).$$ It follows from the homogenity and the monotonicity of the $p$-capacity that 
$$
   C_p(K)\leq (1+\varepsilon)^{n-p}\cdot C_p(int(K)).
$$ 
Hence $C_p(int(K))= C_p(K)$  for all $K\in \cK_0$ by letting $\varepsilon\rightarrow 0^+$.  Please see  \cite{EvansGariepy} for more properties.

Following the convention in the literature of $p$-capacity, in later context we will work on convex domains containing the origin, i.e.,  all open subsets $\Omega\subset \Rn$ whose closure $\overline{\Omega}\in \cK_0$. For convenience, we use $\mathcal{C}_0$ to denote the set of  all open convex domains containing the origin. Moreover, geometric notations for $\Omega\in \cC_0$, such as the support function and the surface area measure,  are considered to be the ones for its closure, for instance, 
$$
  h_{\Omega}(u)=\sup_{x\in \Omega} \langle x, u\rangle=h_{\overline{\Omega}}(u)\ \ \ \mbox{for} \ \ u\in \sphere.
$$  

There exists the $p$-capacitary measure of  $\Omega\in \cC_0$,  denoted by $\mu_p(\Omega, \cdot)$,  on $S^{n-1}$ such that for any Borel set $\Sigma \subset S^{n-1}$ (see e.g., \cite{Lewis1977, LewisNystrom2007, LewisNystrom2008a}), 
\begin{equation}\label{Def-p-capacity-surface area}
   \mu_p(\Omega, \Sigma )=\int_{\mathrm{g}^{-1}(\Sigma )}\|\nabla U_{\Omega}\|^p\, d\mathcal{H}^{n-1},
\end{equation}
where $\mathrm{g}^{-1}:S^{n-1}\to \partial \Omega$ is the inverse Gauss map (i.e., $\mathrm{g}^{-1}(u)$ contains all points $x\in \partial \Omega$ such that $u$ is an unit outer normal vector of $x$) and $U_{\Omega} $ is the $p$-equilibrium potential of $\Omega$. Note that $U_\Omega$ is  the unique solution to the boundary value problem of the following $p$-Laplace equation
$$
   \begin{cases}
     \mbox{div}\left(\|\nabla U \|^{p-2}\nabla U\right)=0 &\mbox{in }\mathbb{R}^n\setminus \overline{\Omega},\\
     U=1 &\mbox{on }\partial \Omega,\\
     \lim_{\|x\|\to \infty}U(x)=0.
   \end{cases}
$$ 
With the help of the $p$-capacitary measure, the Poincar\'e $p$-capacity formula  \cite{CNSXYZ}  gives  $$C_p(\Omega)=\frac{p-1}{n-p}\int_{S^{n-1}} h_\Omega(u)\,d\mu_p(\Omega, u).$$ Lemma 4.1 in \cite{CNSXYZ} asserts that $\mu_p(\Omega_m,\cdot)$ converges to $\mu_p(\Omega,\cdot)$ weakly  on $S^{n-1}$ and hence $C_p(\Omega_m)$ converges to $C_p(\Omega)$, if $\Omega_m$ converges to $\Omega$ in the Hausdorff metric.

The beautiful Hadamard variational formula for $C_p(\cdot)$ was provided in \cite{CNSXYZ}:  for two convex domains $\Omega, \Omega_1\in \cC_0$, one has   
\begin{equation}\label{H-variational-111}
   \frac{1}{n-p}\cdot  \frac{\,d C_p(\Omega+\varepsilon \Omega_1)}{\,d\varepsilon}\bigg|_{\varepsilon=0}=\frac{p-1}{n-p}\int_{S^{n-1}}h_{\Omega_1}(u)\,d\mu_p(\Omega,u)=: C_p(\Omega, \Omega_1), 
\end{equation}  
where $C_p(\Omega, \Omega_1)$ is called the mixed $p$-capacity of $\Omega$ and $\Omega_1$.  By (\ref{class-B-M-capacity-20160615}) and (\ref{H-variational-111}), one gets the $p$-capacitary Minkowski inequality  
\begin{equation}\label{Minkowski-classical-1}
    C_p(\Omega, \Omega_1)^{n-p}\geq C_p(\Omega)^{n-p-1}C_p(\Omega_1),
\end{equation}  
with equality if and only if $\Omega$ and $\Omega_1$ are homothetic \cite{CNSXYZ}. It is also well known that the centroid of $\mu_p(\Omega, \cdot)$ is $o$, that is, 
$$
   \int_{\sphere} u \, d\mu_p(\Omega, u)=o.
$$ 
Moreover, the support of 
$\mu_p(\Omega, \cdot)$ is not contained in any closed hemisphere, i.e.,  there exists a constant $c>0$  such that  
\begin{equation} \label{non-concentration--1}  
    \int_{\sphere} \langle \theta, u\rangle_+ \,d\mu_p(\Omega, u)>c\ \ \ \ \ \mathrm{for\ \ each}\ \ \theta\in \sphere,
\end{equation} 
where  $a_+$ denotes $\max\{a, 0\}$ for all $a\in \R$.

For $\Omega\in \cC_0$ and a Borel set $\Sigma\subset\sphere$, let $$\mu_p(\Omega,\Sigma)=\int_{\Sigma} \, d\mu_p(\Omega, u) \ \ \ \mathrm{and} \ \ \  S(\overline{\Omega},\Sigma)=\int_{\Sigma} \, dS(\overline{\Omega}, u).$$ The following lemma is needed to solve the $p$-capacitary Orlicz-Minkowski problems. See \cite{AGHLV} for a more quantitative argument. 
\bl\label{relation between capacity measure and surface area measure}
Let $\Omega\in \cC_0$ be a convex domain and  $1<p<n$. For a Borel set  $\Sigma\subset \sphere$,   $\mu_p(\Omega,\Sigma)$ and $S(\overline{\Omega},\Sigma)$ either are both strictly positive or are both equal to $0$.  
\el
\begin{proof} 
Let us recall Lemma 2.18 in \cite{CNSXYZ}:  if $\Omega$ is a convex domain such that $\Omega$ is contained in the ball $B(o, R)$ (centered at the origin with radius $R$), there exists a constant $\gamma=\gamma(n,p,R)\in (0, 1]$ such that $\|\nabla U_{\Omega}\| \geq \gamma$ almost everywhere on $\partial \Omega$ (with respect to  $\cH^{n-1}$). This together with formulas (\ref{Def:surface area----123}) and (\ref{Def-p-capacity-surface area}) yield that, for all Borel set $\Sigma\subset \sphere$,  
\begin{eqnarray*}
   \mu_p(\Omega,\Sigma) = \int_{\mathrm{g}^{-1}(\Sigma)} \| \nabla U_{\Omega}(x)\|^{p}\, d\cH^{n-1}(x) \geq   \gamma^p\cdot S(\overline{\Omega},\Sigma).
\end{eqnarray*}  
Consequently if $\mu_p(\Omega,\Sigma)=0$ then $S(\overline{\Omega},\Sigma)=0$ and if $S(\overline{\Omega},\Sigma)>0$   then  $\mu_p(\Omega,\Sigma)>0$.

On the other hand, assume that $S(\overline{\Omega}, \Sigma)=0$ which imples $\cH^{n-1}(\mathrm{g}^{-1}(\Sigma))=0$. Together with formula (\ref{Def-p-capacity-surface area}) and  the fact that $\| \nabla U_{\Omega}\|^{p} $ is integrable on $\partial \Omega$, one has  \begin{eqnarray*}
    \mu_p(\Omega,\Sigma) = \int_{\mathrm{g}^{-1}(\Sigma)} \| \nabla U_{\Omega}(x)\|^{p}\, d\cH^{n-1}(x)=0.
\end{eqnarray*}  
That  is, if $S(\overline{\Omega}, \Sigma)=0$ then $\mu_p(\Omega,\Sigma)=0$ and  if $\mu_p(\Omega,\Sigma) >0$ then $S(\overline{\Omega}, \Sigma)>0$. 
\end{proof}

For $f\in C^{+}(S^{n-1})$,  denote by $\Omega_f$ the Aleksandrov domain associated to $f$ (i.e., the interior of the Aleksandrov body associated to $f$). For $\Omega\in \cC_0$ and $f\in C^{+}(S^{n-1})$, define the mixed $p$-capacity of $\Omega$ and $f$ by  
\begin{equation}\label{mixed-alex-1} 
   C_p(\Omega, f)=\frac{p-1}{n-p}\int_{S^{n-1}}f(u)\ d\mu_p(\Omega, u).
\end{equation}  
Clearly $C_p(\Omega, h_L)=C_p(\Omega, L)$ and $C_p(\Omega, h_\Omega)=C_p(\Omega)$ for all $\Omega, L \in \cC_0$. Moreover, \begin{equation}  \label{capacity-volume-function}
   C_p(\Omega_f)=C_p(\Omega_f, f) 
\end{equation} 
holds for any $f\in C^{+}(S^{n-1})$. This is an immediate consequence of Lemma \ref{relation between capacity measure and surface area measure} (see also \cite[(5.11)]{CNSXYZ}).

\section{The Orlicz $L_{\phi}$  mixed $p$-capacity and related Orlicz-Minkowski inequality} \label{Section:OM-Inequality}

This section is dedicated to prove the $p$-capacitary Orlicz-Hadamard variational formula and establish the $p$-capacitary Orlicz-Minkowski inequality. Let $\phi: (0, \infty)\rightarrow (0, \infty)$ be a continuous function. We now define the Orlicz $L_{\phi}$ mixed $p$-capacity. The mixed $p$-capacity defined in (\ref{H-variational-111}) is related to $\phi=t$.

\bd 
Let $\Omega, \Omega_1\in \cC_0$ be two convex domains. Define $C_{p, \phi}(\Omega, \Omega_1)$, the Orlicz $L_{\phi}$ mixed $p$-capacity of $\Omega$ and $\Omega_1$, by
\begin{equation} \label{def:mixed:orlicz:capacity} 
   C_{p, \phi}(\Omega, \Omega_1)=\frac{p-1}{n-p}\int_{S^{n-1}} \phi\left(\frac{h_{\Omega_1}(u)}{h_{\Omega}(u)}\right)h_{\Omega}(u)\,d\mu_p(\Omega, u).
\end{equation} 
\ed

When $\Omega$ and $\Omega_1$ are dilates of each other, say $\Omega_1=\lambda \Omega$ for some $\lambda>0$, one has 
\begin{equation} \label{dilate-111} 
   C_{p, \phi}(\Omega, \lambda \Omega)= \phi\left(\lambda \right) C_p(\Omega).
\end{equation}

Let  $\varphi_1$ and $\varphi_2$ be either both in $\mathcal{I}$ or both in $\mathcal{D}$. For  $\varepsilon>0$,  let $g_{\varepsilon}$ be defined as in (\ref{definition-linear-orlicz-extension-1}). That is, for $\Omega, \Omega_1\in \cC_0$ and for $u\in \sphere$, 
\begin{equation*}
   \varphi_1\left(\frac{h_{\Omega}(u)}{g_{\varepsilon}(u)}\right)+\varepsilon\varphi_2
   \left(\frac{h_{\Omega_1}(u)}{g_{\varepsilon} (u)}\right)=1.
\end{equation*}  
Clearly  $g_{\varepsilon}\in C^+(\sphere)$. Denote by $\Omega_{\varepsilon}\in \cC_0$ the Aleksandrov domain associated to $g_{\varepsilon}$.

The following lemma for convex domains is identical to Lemma 5.1  in \cite{ZHY2016}.  
\bl \label{convergence 2-0718-1} 
Let $\Omega, \Omega_1 \in\mathcal{C}_0$ and $\varphi_1,\varphi_2\in \mathcal{I}$ be such that $(\varphi_1)'_l(1)$ exists and is positive. Then 
\begin{eqnarray}
   {(\varphi_1)'_l(1)} \lim_{\varepsilon\rightarrow
   0^+}\frac{g_\varepsilon(u)-h_{\Omega} (u)}{\varepsilon}
   &=&{h_{\Omega}(u)}\cdot \varphi_2\left(\frac{h_{\Omega_1}(u)}{h_{\Omega}(u)}\right) \ \ \ \mathrm{uniformly\ on}\ \  S^{n-1}. \label{convergence of f-e}
\end{eqnarray}
For $\varphi_1,\varphi_2\in \mathcal{D}$,  (\ref{convergence of f-e}) holds with $(\varphi_1)'_l(1)$ replaced by $(\varphi_1)'_r(1)$. 
\el

From Lemma \ref{convergence 2-0718-1}, one sees that $g_{\varepsilon}$ converges to $h_\Omega$ uniformly on $\sphere$.  According to  the Aleksandrov convergence lemma,  $\Omega_{\varepsilon}$ converges to $\Omega$ in the Hausdorff metric. We are now ready to establish the geometric interpretation for  the Orlicz $L_{\phi}$ mixed $p$-capacity.  Formula (\ref{H-variational-111}) is the special case when $\varphi_1=\varphi_2=t$.  

\bt \label{c-interpretation-111} 
Let $\Omega, \Omega_1\in \cC_0$ be two convex domains. Suppose  $\varphi_1, \varphi_2\in \mathcal{I}$ such that $ (\varphi_1)'_l(1)$  exists and is nonzero.  Then
\begin{eqnarray*}  
   C_{p, \varphi_2}(\Omega, \Omega_1)=
   \frac{{(\varphi_1)'_l(1)}}{n-p}\cdot \lim_{\varepsilon\to 0^+}\frac{C_p(\Omega_{\varepsilon})-C_p(\Omega)}{\varepsilon}.
\end{eqnarray*} 
With  $(\varphi_1)'_l(1)$ replaced by $(\varphi_1)'_r(1)$ if $(\varphi_1)'_r(1)$  exists and is nonzero, one gets the analogous result  for $\varphi_1, \varphi_2\in \mathcal{D}$. \et

\begin{proof} 
The proof of this theorem is similar to analogous results in \cite{CNSXYZ, GHW2014, HLYZ, XJL, ZHY2016}. A brief proof is included  here for completeness. As $\Omega_\varepsilon\rightarrow \Omega$ in the Hausdorff metric,  $\mu_p(\Omega_\varepsilon,\cdot)\rightarrow \mu_p(\Omega,\cdot)$ weakly  on $S^{n-1}$ due to  Lemma 4.1 in \cite{CNSXYZ}. Moreover, if $h_{\varepsilon}\rightarrow h$ uniformly on $\sphere$, then 
$$
   \lim_{\varepsilon\rightarrow 0^+} \int_{\sphere} h_{\varepsilon}(u) \,d\mu_p(\Omega_\varepsilon, u)=\int_{\sphere} h(u) \,d\mu_p(\Omega, u).
$$
In particular, it follows from (\ref{mixed-alex-1}) and Lemma \ref{convergence 2-0718-1}  that 
\begin{eqnarray*}
   (\varphi_1)_{l}^{'}(1) \cdot \lim_{\varepsilon\rightarrow 0^+}\frac{C_p(\Omega_\varepsilon, g_\varepsilon)-C_p(\Omega_\varepsilon, h_{\Omega})}{\varepsilon}
   &=&(\varphi_1)_{l}^{'}(1) \cdot 
   \lim_{\varepsilon\rightarrow 0^+}\frac{p-1}{n-p}\int_{S^{n-1}}\frac{g_\varepsilon(u)-h_{\Omega}(u)}{\varepsilon}d\mu_{p}(\Omega_\varepsilon,u)\\
   &=& \frac{p-1}{n-p}\int_{S^{n-1}} h_{\Omega}(u)  \varphi_2 \left(\frac{h_{\Omega_1}(u)}{h_{\Omega}(u)}\right)d\mu_{p}(\Omega, u)\\
   &=& C_{p,\varphi_2}(\Omega, \Omega_1). 
\end{eqnarray*} 
Inequality  (\ref{Minkowski-classical-1}), formula (\ref{capacity-volume-function}), and the continuity of $p$-capacity yield that
\begin{eqnarray*}
   C_{p,\varphi_2}(\Omega, \Omega_1) &=&(\varphi_1)_{l}^{'}(1)\cdot \liminf_{\varepsilon\rightarrow 0^+}\frac{C_p(\Omega_\varepsilon)-C_p(\Omega_\varepsilon, \Omega)}{\varepsilon} \\ &\leq&(\varphi_1)_{l}^{'}(1)\cdot
   \liminf_{\varepsilon\rightarrow 0^+} \bigg[C_p(\Omega_\varepsilon)^{\frac{n-p-1}{n-p}}\cdot \frac{C_p(\Omega_\varepsilon)^{\frac{1}{n-p}}-C_p(\Omega)^{\frac{1}{n-p}}}{\varepsilon}\bigg] \\
   &=& (\varphi_1)_{l}^{'}(1)\cdot C_p(\Omega)^{\frac{n-p-1}{n-p}}\cdot \liminf_{\varepsilon\rightarrow 0^+}\frac{C_p(\Omega_\varepsilon)^{\frac{1}{n-p}}-C_p(\Omega)^{\frac{1}{n-p}}}{\varepsilon}.
\end{eqnarray*}  
Similarly,  as  $h_{\Omega_\varepsilon} \leq g_\varepsilon$ and $C_p(\Omega)=C_p(\Omega, h_{\Omega})$, one has  
\begin{eqnarray*}
   C_{p,\varphi_2}(\Omega, \Omega_1) &=& (\varphi_1)_{l}^{'}(1)\cdot\lim_{\varepsilon\rightarrow 0^+}\frac{p-1}{n-p}\int_{S^{n-1}}\frac{g_\varepsilon(u)-h_{\Omega}(u)}{\varepsilon}d\mu_{p}(\Omega, u)\\ &\geq& (\varphi_1)_{l}^{'}(1)\cdot\limsup_{\varepsilon\rightarrow0^+}\frac{C_p(\Omega, \Omega_\varepsilon)-C_p(\Omega)}{\varepsilon}\\ &\geq&   (\varphi_1)_{l}^{'}(1)\cdot C_p(\Omega)^{\frac{n-p-1}{n-p}}\cdot\limsup_{\varepsilon\rightarrow0^+}\frac{C_p(\Omega_\varepsilon)^{\frac{1}{n-p}}-C_p(\Omega)^{\frac{1}{n-p}}}{\varepsilon}.
\end{eqnarray*} 
This concludes that   
\begin{eqnarray*}
   C_{p,\varphi_2}(\Omega, \Omega_1) &=&  (\varphi_1)_{l}^{'}(1)\cdot C_p(\Omega)^{\frac{n-p-1}{n-p}}\cdot\lim_{\varepsilon\rightarrow 0^+}\frac{C_p(\Omega_\varepsilon)^{\frac{1}{n-p}}-C_p(\Omega)^{\frac{1}{n-p}}}{\varepsilon} \\ &=&
   \frac{{(\varphi_1)'_l(1)}}{n-p}\cdot \lim_{\varepsilon\to 0^+}\frac{C_p(\Omega_{\varepsilon})-C_p(\Omega)}{\varepsilon},
\end{eqnarray*} 
where the second equality follows from a standard argument by the chain rule.   
\end{proof}

Let  $p\in (1, n)$ and $q\neq 0$ be real numbers. For $\Omega, \Omega_1\in \cC_0$, define $C_{p, q}(\Omega, \Omega_1)$, the $L_{q}$ mixed $p$-capacity of $\Omega$ and $\Omega_1$, by
\begin{equation} \label{def:mixed:orlicz:capacity-qq} 
   C_{p, q}(\Omega, \Omega_1)=\frac{p-1}{n-p}\int_{S^{n-1}}   \big[h_{\Omega_1}(u)\big]^q  \,d\mu_{p, q} (\Omega, u),
\end{equation} 
where $\mu_{p, q}(\Omega, \cdot)$ denotes the $L_q$ $p$-capacitary measure of $\Omega$: 
\begin{equation*} 
   \,d\mu_{p, q}(\Omega, \cdot)=h_{\Omega}^{1-q} \,d\mu_p(\Omega, \cdot).  
\end{equation*}    
For $\varepsilon>0$, let $h_{q, \varepsilon}=\big[h_{\Omega}^q+\varepsilon h_{\Omega_1}^q\big]^{1/q}$  and $\Omega_{h_{q, \varepsilon}}$ be the  Aleksandrov domain associated to $h_{q, \varepsilon}$.  By letting $\varphi_1=\varphi_2=t^q$ for $q\neq 0$ in Theorem \ref{c-interpretation-111}, one gets the geometric interpretation for  $C_{p, q}(\cdot, \cdot)$. 
\bc \label{c-q-interpretation-111} 
Let $\Omega, \Omega_1\in \cC_0$  and $p\in (1, n)$. For all $0\neq q\in \R$, one has 
\begin{eqnarray*}  
 	C_{p, q}(\Omega, \Omega_1)=
    \frac{{q}}{n-p}\cdot \lim_{\varepsilon\to 0^+}\frac{C_p(\Omega_{h_{q, \varepsilon}})-C_p(\Omega)}{\varepsilon}.
\end{eqnarray*}  
\ec

Regarding the Orlicz $L_{\phi}$ mixed $p$-capacity, one has the following $p$-capacitary Orlicz-Minkowski inequality. When $\phi=t$, one recovers the $p$-capacitary Minkowski inequality (\ref{Minkowski-classical-1}).

\bt \label{capacitary-Minkowski-2-2-1-1} 
Let $\Omega, \Omega_1\in \cC_0$ and $p\in (1, n)$. Suppose that   $\phi: [0, \infty)\rightarrow [0, \infty)$ is increasing and convex.  Then
\begin{eqnarray*}   
	C_{p, \phi}(\Omega, \Omega_1)  \geq C_p(\Omega) \cdot  \phi\bigg(\bigg(\frac{C_p(\Omega_1)}{C_p(\Omega)}\bigg)^{\frac{1}{n-p}}\bigg).
\end{eqnarray*}
If in addition $\phi$ is strictly convex,   equality holds if and only if $\Omega$ and $\Omega_1$ are dilates of each other. 
\et 

\begin{proof} 
It follows from Jensen's inequality (see \cite{ghwy2015}), $C_{p}(\Omega)>0$ and the convexity of $\phi$ that
\begin{eqnarray}  
C_{p, \phi}(\Omega, \Omega_1) &=& \frac{p-1}{n-p} \int_{S^{n-1}} \phi\left(\frac{h_{\Omega_1}(u)}{h_{\Omega}(u)}\right)h_{ \Omega}(u)\,d\mu_p(\Omega, u) \nonumber \\& \geq & C_p(\Omega) \cdot \phi\left(\int_{S^{n-1}} \frac{p-1}{n-p }\cdot \frac{h_{\Omega_1}(u)} {C_p(\Omega)}\,d\mu_p(\Omega, u)\right) \nonumber \\ &=&C_p(\Omega) \cdot \phi\left(\frac{C_{p, 1}(\Omega, \Omega_1)} {C_p(\Omega)}\right) \nonumber
\\ &\geq & C_p(\Omega) \cdot  \phi\bigg(\bigg(\frac{C_p(\Omega_1)}{C_p(\Omega)}\bigg)^{\frac{1}{n-p}}\bigg)\label{Minkowski-Orlicz--1-1---1-1-1-1-1}
\end{eqnarray} 
where the last inequality follows from (\ref{Minkowski-classical-1}) and the fact that $\phi$ is increasing.

From (\ref{homogeneous-111111111111}) and (\ref{dilate-111}), if $\Omega$ and $\Omega_1$ are dilates of each other, then clearly \begin{eqnarray*}   
	C_{p, \phi}(\Omega, \Omega_1)  = C_p(\Omega) \cdot  \phi\bigg(\bigg(\frac{C_p(\Omega_1)}{C_p(\Omega)}\bigg)^{\frac{1}{n-p}}\bigg).
\end{eqnarray*} 
On the other hand, if  $\phi$ is strictly convex,  equality holds in (\ref{Minkowski-Orlicz--1-1---1-1-1-1-1}) only if equalities hold in
both the first and the second inequalities of  (\ref{Minkowski-Orlicz--1-1---1-1-1-1-1}). For the second one, $\Omega$ and $\Omega_1$ are homothetic to each other. That is, there exists $r>0$ and $x\in \Rn$, such that $\Omega_1=r\Omega+x$ and hence for all $u\in \sphere$,  
$$
   h_{\Omega_1}(u)=r \cdot h_{\Omega}(u)+\langle x, u\rangle. 
$$
As $\phi$ is strictly convex, the characterization of equality in Jensen's inequality implies that 
$$
   \frac{h_{\Omega_1}(v)}{h_{\Omega}(v)}= \int_{S^{n-1}} \frac{p-1}{n-p }\cdot \frac{h_{\Omega_1}(u)} {C_p(\Omega)}\,d\mu_p(\Omega, u)
$$ 
for $\mu_p(\Omega, \cdot)$-almost all $v\in \sphere$. This together with the fact that  $\mu_p(\Omega, \cdot)$ has its centroid at the origin yield  $\langle x, v\rangle=0$  for $\mu_p(\Omega, \cdot)$-almost all $v\in \sphere$. As the support of $\mu_p(\Omega, \cdot)$ is not  contained in any closed hemisphere, one has $x=0$. That is, $\Omega$ and $\Omega_1$ are dilates of each other.   \end{proof}

An application of the above $p$-capacitary Orlicz-Minkowski inequality is stated below.
\bt  \label{uniqueness-mixed volume-1}  
Let $\phi\in \Phi_1$ be strictly increasing and strictly convex.   Assume that $\Omega, \widetilde{\Omega}\in \cC_0$ are two convex domains. Then $\Omega=\widetilde{\Omega}$  if the following equality holds for all $\Omega_1\in \cC_0$:   
\begin{equation}\label{condition---1---1}
   \frac{C_{p, \phi}(\Omega, \Omega_1)}{C_p(\Omega)}=\frac{C_{p, \phi}(\widetilde{\Omega}, \Omega_1)}{C_p(\widetilde{\Omega})}.
\end{equation} 
Moreover, $\Omega=\widetilde{\Omega}$ also holds if, for any $\Omega_1\in \cC_0$,  
\begin{equation}\label{condition---2---2} 
   C_{p, \phi}(\Omega_1, \Omega) = C_{p, \phi}(\Omega_1, \widetilde{\Omega}).
\end{equation}  
\et 
\begin{proof} 
It follows from equality (\ref{condition---1---1}) and the $p$-capacitary Orlicz-Minkowski inequality that  
\begin{equation} \label{equa-1-1--1}
	1=\frac{C_{p, \phi}(\Omega, \Omega)}{C_p(\Omega)}=\frac{C_{p, \phi}(\widetilde{\Omega}, \Omega)}{C_p(\widetilde{\Omega})}\geq \phi\bigg(\bigg(\frac{C_p(\Omega)}{C_p(\widetilde{\Omega})}\bigg)^{\frac{1}{n-p}}\bigg). 
\end{equation} 
The fact that $\phi$ is strictly increasing with $\phi(1)=1$ and $n-p>0$ yield  $C_p(\widetilde{\Omega})\geq C_p(\Omega)$. Similarly,  $C_p(\widetilde{\Omega})\leq C_p(\Omega)$ and then $C_p(\widetilde{\Omega}) = C_p(\Omega)$. Hence,  equality holds in inequality (\ref{equa-1-1--1}).  This can happen only if $\Omega$ and $\widetilde{\Omega}$ are dilates of each other, due to Theorem \ref{capacitary-Minkowski-2-2-1-1} and the fact that $\phi$ is strictly convex. Combining with the above proved fact $C_p(\widetilde{\Omega})= C_p(\Omega)$, one gets  $\Omega=\widetilde{\Omega}$.  

Follows along the same lines, $\Omega=\widetilde{\Omega}$ if equality (\ref{condition---2---2}) holds for any $\Omega_1\in \cC_0$.
\end{proof}

Note that $\phi=t^q$ for $q>1$ is a strictly convex and strictly increasing function. Theorem  \ref{capacitary-Minkowski-2-2-1-1} yields the $p$-capacitary $L_q$ Minkowski inequality:  for  $\Omega, \Omega_1\in \cC_0$, one has 
\begin{eqnarray*}  
	C_{p, q}(\Omega, \Omega_1)  \geq \big[C_p(\Omega)\big]^{\frac{n-p-q}{n-p}} \cdot  \big[C_p(\Omega_1)\big]^{\frac{q}{n-p}}
\end{eqnarray*} 
with equality if and only if $\Omega$ and $\Omega_1$ are dilates of each other.  
\bc \label{uniqueness-measure-1} 
Let $p\in (0, n)$ and  $q>1$.  If $\Omega, \widetilde{\Omega}\in \cC_0$ are such that 
$$
   \mu_{p, q}(\Omega, \cdot)=\mu_{p, q}(\widetilde{\Omega}, \cdot),
$$   
then $\Omega=\widetilde{\Omega}$ if $q\neq n-p$, and $\Omega$ is dilate of $\widetilde{\Omega}$ if $q=n-p$. 
\ec 
\begin{proof} 
Firstly let $q>1$ and $q\neq n-p$. As $\mu_{p, q}(\Omega, \cdot)=\mu_{p, q}(\widetilde{\Omega}, \cdot)$, it follows form (\ref{def:mixed:orlicz:capacity-qq}) that, for all $\Omega_1\in \cC_0$,  \begin{eqnarray}\label{p-q-equality-1030} 
	C_{p, q}(\Omega, \Omega_1)=C_{p, q}(\widetilde{\Omega}, \Omega_1). 
\end{eqnarray}   
By letting $\Omega_1=\widetilde{\Omega}$, one has, 
\begin{eqnarray*}   
	C_{p, q}(\Omega, \widetilde{\Omega}) =C_p(\widetilde{\Omega}) \geq \big[C_p(\Omega)\big]^{\frac{n-p-q}{n-p}} \cdot  \big[C_p(\widetilde{\Omega})\big]^{\frac{q}{n-p}}.
\end{eqnarray*} 
This yields $C_p(\Omega)\geq C_p(\widetilde{\Omega}) $ if $q>n-p$ and  $C_p(\Omega)\leq C_p(\widetilde{\Omega}) $ if $q<n-p$.  
Similarly, by letting $\Omega_1=\Omega$, one has $C_p(\Omega)\leq C_p(\widetilde{\Omega}) $ if $q>n-p$ and  $C_p(\Omega)\geq C_p(\widetilde{\Omega}) $ if $q<n-p$. In any cases, $C_p(\Omega)= C_p(\widetilde{\Omega})$. Together with (\ref{p-q-equality-1030}), Theorem  \ref{uniqueness-mixed volume-1} yields the desired argument $\Omega=\widetilde{\Omega}$. 

Now assume that $q=n-p>1$. Then  (\ref{p-q-equality-1030}) yields \begin{eqnarray*}   
	C_{p, q}(\Omega, \widetilde{\Omega}) =C_p(\widetilde{\Omega}) \geq \big[C_p(\Omega)\big]^{\frac{n-p-q}{n-p}} \cdot  \big[C_p(\widetilde{\Omega})\big]^{\frac{q}{n-p}}=C_p(\widetilde{\Omega}).
\end{eqnarray*} 
It follows from Theorem \ref{capacitary-Minkowski-2-2-1-1} that $\Omega$ and $\widetilde{\Omega}$ are dilates of each other. 
\end{proof} 
 
It is worth to mention that $C_{p, \phi}(\cdot, \cdot)$ is not homogeneous if $\phi$ is not a homogeneous function; this can be seen from formula (\ref{dilate-111}). When $\phi\in \mathcal{I}$, we can define $\widehat{C}_{p, \phi}(\Omega, \Omega_1)$, the homogeneous Orlicz $L_{\phi}$ mixed $p$-capacity of $\Omega, \Omega_1\in \cC_0$,  by 
$$
   \widehat{C}_{p, \phi}(\Omega, \Omega_1)=\inf\bigg\{\eta>0:   \frac{p-1}{n-p}\int_{S^{n-1}} \phi\left(\frac{h_{\Omega_1}(u)}{\eta\cdot h_{\Omega}(u)}\right)h_{\Omega}(u)\,d\mu_p(\Omega, u) \leq C_{p}(\Omega)\bigg\},
$$ 
while $\widehat{C}_{p, \phi}(\Omega, \Omega_1)$ for $\phi\in \mathcal{D}$ is defined as above with ``$\leq$" replaced by ``$\geq$". If $\phi=t^q$ for $q\neq 0$,  
$$
   \widehat{C}_{p, \phi}(\Omega, \Omega_1)=\bigg(\frac{C_{p, q}(\Omega, \Omega_1)}{C_p(\Omega)}\bigg)^{1/q}. 
$$ 
For all $\eta>0$ and for $\phi\in \mathcal{I}$, let 
$$
   g(\eta)=\frac{p-1}{n-p} \int_{S^{n-1}} \phi\left(\frac{h_{\Omega_1}(u)}{\eta\cdot h_{\Omega}(u)}\right) h_{\Omega}(u)\,d\mu_p(\Omega, u).
$$ 
The fact that $\phi$ is monotone increasing yields 
\begin{eqnarray*}  
	\phi\left(\frac{\min_{u\in \sphere} h_{\Omega_1}(u)}{\eta\cdot \max_{u\in\sphere} h_{\Omega}(u)}\right) \leq \frac{g(\eta)}{C_p(\Omega)} \leq  \phi\left(\frac{\max_{u\in \sphere} h_{\Omega_1}(u)}{\eta\cdot \min_{u\in\sphere} h_{\Omega}(u)}\right).
\end{eqnarray*} 
Hence $\lim_{\eta\rightarrow 0^+}g(\eta)=\infty$ and $\lim_{\eta\rightarrow \infty}g(\eta)=0.$ It is also easily checked that $g$ is strictly decreasing. This concludes that if $\phi\in \mathcal{I}$, 
\begin{equation} \label{equation-equivalent-homogeneous-mixed-capacity}  
   \frac{p-1}{n-p}\int_{S^{n-1}} \phi\left(\frac{h_{\Omega_1}(u)}{\widehat{C}_{p, \phi}(\Omega, \Omega_1) \cdot h_{\Omega}(u)}\right)h_{\Omega}(u)\,d\mu_p(\Omega, u) = C_{p}(\Omega). 
\end{equation}   
Following along the same lines, formula (\ref{equation-equivalent-homogeneous-mixed-capacity}) also holds for $\phi\in \mathcal{D}$. 
 
The $p$-capacitary Orlicz-Minkowski inequality for $\widehat{C}_{p, \phi}(\cdot, \cdot)$ is stated in the following result. 
\bc 
Let $\phi\in \mathcal{I}$ be convex. For all $\Omega, \Omega_1\in \cC_0$, one has,  
\begin{equation} \label{homogeneous-Minkowski--1-1--1-1-1} 
   \widehat{C}_{p, \phi}(\Omega, \Omega_1)\geq \bigg(\frac{C_p(\Omega_1)}{C_p(\Omega)}\bigg)^{\frac{1}{n-p}}.
\end{equation}  
If in addition $\phi$ is strictly convex, equality holds if and only if $\Omega$ and $\Omega_1$ are dilates of each other.  
\ec 
\begin{proof}  
It follows from formula (\ref{equation-equivalent-homogeneous-mixed-capacity})  and Jensen's inequality that 
\begin{eqnarray*}   
 	1&=& \int_{S^{n-1}} \phi\bigg(\frac{h_{\Omega_1}(u)}{\widehat{C}_{p, \phi}(\Omega, \Omega_1)  \cdot h_{\Omega}(u)}\bigg) \cdot \frac{p-1}{n-p}  \cdot \frac{h_{\Omega}(u)}{C_{p}(\Omega)}\,d\mu_p(\Omega, u)  \\ &\geq& \phi \bigg(\int_{S^{n-1}} \frac{h_{\Omega_1}(u)}{\widehat{C}_{p, \phi}(\Omega, \Omega_1) } \cdot \frac{p-1}{n-p}  \cdot \frac{1}{C_{p}(\Omega)}\,d\mu_p(\Omega, u) \bigg)\\ &=& \phi\bigg(\frac{C_p(\Omega, \Omega_1)}{\widehat{C}_{p, \phi}(\Omega, \Omega_1) \cdot C_{p}(\Omega)}\bigg).  
\end{eqnarray*} 
As $\phi(1)=1$ and $\phi$ is monotone increasing, one has $$\widehat{C}_{p, \phi}(\Omega, \Omega_1) \geq \frac{C_p(\Omega, \Omega_1)}{ C_{p}(\Omega)}\geq \bigg(\frac{C_p(\Omega_1)}{C_p(\Omega)}\bigg)^{\frac{1}{n-p}},$$ where the second inequality follows from (\ref{Minkowski-classical-1}). 
 
 It is easily checked that equality holds in (\ref{homogeneous-Minkowski--1-1--1-1-1}) if $\Omega_1$ is dilate of $\Omega$. Now assume that in addition $\phi$ is strictly convex and equality holds in (\ref{homogeneous-Minkowski--1-1--1-1-1}). Then equality must hold in (\ref{Minkowski-classical-1}) and hence $\Omega$ is homothetic to $\Omega_1$. Following along the same lines in the proof of Theorem \ref{capacitary-Minkowski-2-2-1-1}, one obtains that $\Omega$ is dilate of $\Omega_1$. 
 \end{proof}

\section{The $p$-capacitary Orlicz-Brunn-Minkowski inequality}\label{Section:OBM}

This section aims to establish the $p$-capacitary Orlicz-Brunn-Minkowski inequality (i.e., Theorem \ref{C-O-B-M-inequality-1}). We also show that the $p$-capacitary Orlicz-Brunn-Minkowski inequality is equivalent to the $p$-capacitary  Orlicz-Minkowski inequality (i.e., Theorem \ref{capacitary-Minkowski-2-2-1-1}) in some sense. Let $m\geq 2$.  Recall that the support function of $+_{\varphi}(\Omega_1,\dots, \Omega_m)$ satisfies the following equation: for any $u\in \sphere$,   \begin{equation}\label{Orldef:1--1--1}
   \varphi\bigg(\frac{h_{\Omega_1}(u)}{h_{+_{\varphi}(\Omega_1,\dots, \Omega_m)}(u)},\dots, \frac{h_{\Omega_m}(u)}{h_{+_{\varphi}(\Omega_1,\dots, \Omega_m)}(u)}\bigg) = 1.
\end{equation}

\bt \label{C-O-B-M-inequality-1}
Suppose that  $\Omega_1, \cdots, \Omega_m \in \cC_0$ are convex domains. For all  $\varphi\in \Phi_m$, one has  
\begin{equation} \label{Orlicz-B-M-inequality---11111--111}
   1\geq \varphi\bigg(\left(\frac{C_p(\Omega_1)}{C_p(+_\varphi (\Omega_1, \cdots, \Omega_m))}\right)^{\frac{1}{n-p}}, \cdots, \left(\frac{C_p(\Omega_m)}{C_p(+_\varphi (\Omega_1, \cdots, \Omega_m))}\right)^{\frac{1}{n-p}}\bigg).
\end{equation}
If in addition $\varphi$ is strictly convex, equality holds if and only if  $\Omega_i$ are dilates of $\Omega_1$   for all $i=2, 3, \cdots, m$.
\et

\begin{proof} 
Let $\varphi\in \Phi_m$ and  $\Omega_1, \cdots, \Omega_m \in \cC_0$. Recall that $\Omega_1\subset +_{\varphi}(\Omega_1,\dots, \Omega_m)$ (see (\ref{compare-Orliczdef})). The fact that the $p$-capacity is monotone increasing yields 
$$
   C_p(+_{\varphi}(\Omega_1,\dots, \Omega_m))\geq C_p(\Omega_1)>0.
$$ 
Define a probability measure on $S^{n-1}$ by 
$$
d\omega_{p, \varphi} (u)=\frac{p-1}{n-p}\cdot\frac{1}{C_p(+_{\varphi}(\Omega_1,\dots, \Omega_m))}\cdot h_{ +_{\varphi}(\Omega_1,\dots, \Omega_m)}(u) \,d\mu_p(+_{\varphi}(\Omega_1,\dots, \Omega_m), u).
$$ 
It follows from formulas (\ref{def:mixed:orlicz:capacity}) and (\ref{Orldef:1--1--1}),  and Jensen's inequality (see \cite[Proposition 2.2]{ghwy2015}) that \begin{eqnarray*}
   1 &=& \int_{\sphere}\varphi\bigg(\frac{h_{\Omega_1}(u)}{h_{+_{\varphi}(\Omega_1, \cdots, \Omega_m)}(u)}, \cdots, \frac{h_{\Omega_m}(u)}{h_{+_{\varphi}(\Omega_1, \cdots, \Omega_m)}(u)}\bigg) \,d\omega_{p, \varphi} (u)\\ &\geq& \varphi\bigg(\int_{\sphere} \frac{h_{\Omega_1}(u)}{h_{+_{\varphi}(\Omega_1, \cdots, \Omega_m)}(u)}\,d\omega_{p, \varphi} (u),\cdots, \int_{\sphere} \frac{h_{\Omega_m}(u)}{h_{+_{\varphi}(\Omega_1, \cdots, \Omega_m)}(u)}\,d\omega_{p, \varphi} (u) \bigg) \\
   &=&\varphi\left( \frac{C_{p,1}(+_{\varphi}(\Omega_1, \cdots, \Omega_m), \Omega_1)}{C_p(+_{\varphi}(\Omega_1, \cdots, \Omega_m))}, \cdots, \frac{C_{p,1}(+_{\varphi}(\Omega_1, \cdots, \Omega_m), \Omega_m)}{C_p(+_{\varphi}(\Omega_1, \cdots, \Omega_m))}\right)\\
   &\geq& \varphi\bigg(\bigg(\frac{C_p(\Omega_1)}{C_p(+_{\varphi}(\Omega_1, \cdots, \Omega_m))}\bigg)^{\frac{1}{n-p}}, \cdots,  \bigg(\frac{C_p(\Omega_m)}{C_p(+_{\varphi}(\Omega_1, \cdots, \Omega_m))}\bigg)^{\frac{1}{n-p}} \bigg),
\end{eqnarray*} 
where the last inequality follows from inequality (\ref{Minkowski-classical-1}).

Let us now characterize the conditions for equality. In fact, if $\Omega_i$ are dilates of $\Omega_1$ for all $1<i\leq m$, then $+_{\varphi}(\Omega_1, \cdots, \Omega_m)$ is also dilate of $\Omega_1$ and the equality clearly holds.  Now suppose that $\varphi\in \Phi_m$ is strictly convex. Equality must hold for Jensen's inequality and hence there exists a vector $z_0\in \R^m$ (see \cite[Proposition 2.2]{ghwy2015})  such that 
\begin{equation*}  
   \bigg(\frac{h_{\Omega_1}(u)}{h_{+_{\varphi}(\Omega_1, \cdots, \Omega_m)}(u)},\cdots, \frac{h_{\Omega_m}(u)}{h_{+_{\varphi}(\Omega_1, \cdots, \Omega_m)}(u)}\bigg)=z_0
\end{equation*}  
for $\omega_{p, \varphi}$-almost all $u\in S^{n-1}$. Moreover, as $\varphi\in \Phi_m$ is strictly increasing on each component, one must have  
$$ 
   \frac{C_{p,1}(+_{\varphi}(\Omega_1, \cdots, \Omega_m), \Omega_j)}{C_p(+_{\varphi}(\Omega_1, \cdots, \Omega_m))}=\bigg(\frac{C_p(\Omega_j)}{C_p(+_{\varphi}(\Omega_1, \cdots, \Omega_m))}\bigg)^{\frac{1}{n-p}}
$$ 
for all $j=1, 2, \cdots, m$. The characterization of equality for (\ref{Minkowski-classical-1}) yields that  $\Omega_j$ for $j=1, \cdots, m$ are all  homothetic to  $+_{\varphi}(\Omega_1, \cdots, \Omega_m)$. Following the argument similar to that of Theorem \ref{capacitary-Minkowski-2-2-1-1}, one can conclude that $\Omega_i$ for all $j=1, \cdots, m$ are dilates of $+_{\varphi}(\Omega_1, \cdots, \Omega_m)$,  as desired. 
\end{proof}

If $\varphi(x)=\sum_{i=1}^m x_i$ for $x\in [0, \infty)^m$, then $\varphi\in \Phi_m$ and inequality (\ref{Orlicz-B-M-inequality---11111--111}) becomes  the classical $p$-capacitary Brunn-Minkowski inequality (see inequality (\ref{class-B-M-capacity-20160615})): for $\Omega_1, \cdots, \Omega_m\in \cC_0$, one has  
\begin{equation} \label{Brunn-Minkowski-classical-p-capacity} 
   C_p(\Omega_1+\cdots+\Omega_m)^{\frac{1}{n-p}}\geq C_p(\Omega_1)^{\frac{1}{n-p}}+\cdots +C_p(\Omega_m)^{\frac{1}{n-p}}.
\end{equation}  
From the proof of Theorem \ref{C-O-B-M-inequality-1}, one sees that  equality holds if and only if $\Omega_i$ is homothetic to $\Omega_j$ for all $1\leq i<j\leq m$.   When $\varphi(x)=\sum_{i=1}^m x_j^q\in \Phi_m$ for $q>1$, one gets the $p$-capacitary  $L_q$-Brunn-Minkowski inequality:  for $\Omega_1, \cdots, \Omega_m\in \cC_0$, one has 
$$
   C_p(\Omega_1+_q\cdots+_q \Omega_m)^{\frac{q}{n-p}}\geq C_p(\Omega_1)^{\frac{q}{n-p}}+\cdots +C_p(\Omega_m)^{\frac{q}{n-p}}.
$$ 
As $\varphi(x)=\sum_{i=1}^m x_j^q$ for $q>1$ is strictly convex,  equality holds if and only if  $\Omega_i$ is dilate of $\Omega_j$ for all $1\leq i<j\leq m$. This has been proved by Zou and Xiong in \cite{xiongzou} with  a  different approach.

Now let us consider the linear Orlicz addition of $\Omega_1, \cdots, \Omega_m\in \cC_0$. This is related to  \begin{equation}\label{def-linear-varphi}
   \varphi(x)=\alpha_1\varphi_1(x_1)+\cdots +\alpha_m\varphi_m(x_m), \ \ \ x=(x_1,\cdots, x_m)\in (0, \infty)^m,   
\end{equation} 
where $\alpha_j>0$ are constants and $\varphi_j\in \Phi_1$ for all $j=1, \cdots, m$. Clearly $\varphi\in \Phi_m$ and  the $p$-capacitary Orlicz-Brunn-Minkowski inequality in Theorem \ref{C-O-B-M-inequality-1} can be rewritten as the following form.
\bt\label{C-O-B-M-inequality-2}
Let  $\varphi$ be given in (\ref{def-linear-varphi}) with $\alpha_j>0$ constants and $\varphi_j\in \Phi_1$ for $j=1, \cdots, m$. For $\Omega_1, \cdots, \Omega_m \in \cC_0$, one has
\begin{equation} \label{Orlicz-Minkowski-inequality-linear--1} 
   1\geq \sum_{j=1}^m \alpha_j \varphi_j\bigg(\bigg(\frac{C_p(\Omega_j)}{C_p(+_\varphi (\Omega_1, \cdots, \Omega_m))}\bigg)^{\frac{1}{n-p}}\bigg).
\end{equation} 
\et

In fact, inequality (\ref{Orlicz-Minkowski-inequality-linear--1}) is equivalent to, in some sense, the $p$-capacitary Orlicz-Minkowski inequality in Theorem \ref{capacitary-Minkowski-2-2-1-1}.  Let  $m=2$, $\varphi_1, \varphi_2\in \Phi_1$, $\Omega, \widetilde{\Omega} \in \cC_0$, $\alpha_1=1$ and $\alpha_2=\varepsilon>0$. In this case, the linear Orlicz addition of $\Omega$ and $\widetilde{\Omega}$ is denoted by $\Omega+_{\varphi, \varepsilon} \widetilde{\Omega}$, whose support function is given by, for $u\in S^{n-1}$, $$\varphi_1 \bigg(\frac{h_{\Omega}(u)}{h_{\Omega+_{\varphi, \varepsilon}\widetilde{\Omega}}(u)}\bigg)+\varepsilon \varphi_2\bigg(\frac{h_{\widetilde{\Omega}}(u)}{h_{\Omega +_{\varphi, \varepsilon} \widetilde{\Omega}}(u)}\bigg)=1. $$ The $p$-capacitary Orlicz-Brunn-Minkowski inequality in Theorem \ref{C-O-B-M-inequality-2} becomes 
\begin{equation*} 
   1\geq \varphi_1\bigg(\bigg(\frac{C_p(\Omega)}{C_p(\Omega+_{\varphi, \varepsilon} \widetilde{\Omega})}\bigg)^{\frac{1}{n-p}}\bigg)+\varepsilon \varphi_2\bigg (\bigg(\frac{C_p(\widetilde{\Omega})}{C_p(\Omega+_{\varphi, \varepsilon} \widetilde{\Omega})}\bigg)^{\frac{1}{n-p}}\bigg),
\end{equation*} 
for all $\varepsilon>0$. It is equivalent to
\begin{equation}\label{C-O-B-M-inequality-2-111}
   1-\varphi_1^{-1}\bigg(1-\varepsilon \varphi_2\bigg(\bigg(\frac{C_p(\widetilde{\Omega})}{C_p(\Omega+_{\varphi, \varepsilon} \widetilde{\Omega})}\bigg)^{\frac{1}{n-p}}\bigg)\bigg)
   \leq 1- \bigg( \frac{C_p(\Omega)}  {C_p(\Omega+_{\varphi, \varepsilon} \widetilde{\Omega})}\bigg)^{\frac{1}{n-p}}.
\end{equation}  
For convenience, let $z(\varepsilon)$ be 
$$
   z(\varepsilon)=\varphi_1^{-1}\bigg(1-\varepsilon \varphi_2\bigg(\bigg(\frac{C_p(\widetilde{\Omega})}{C_p(\Omega+_{\varphi, \varepsilon} \widetilde{\Omega})}\bigg)^{\frac{1}{n-p}}\bigg)\bigg).
$$ 
Then $z(\varepsilon)\rightarrow 1^-$ as $\varepsilon\rightarrow 0^+$ and \begin{eqnarray*}    
 	\lim_{\varepsilon\rightarrow 0^+} \frac{1-z(\varepsilon)}{\varepsilon}  &=& \lim_{\varepsilon\rightarrow 0^+} \frac{1-z(\varepsilon)}{1-\varphi_1(z(\varepsilon))}  \cdot \lim_{\varepsilon\rightarrow 0^+}\varphi_2\bigg(\bigg(\frac{C_p(\widetilde{\Omega})}{C_p(\Omega+_{\varphi, \varepsilon} \widetilde{\Omega})}\bigg)^{\frac{1}{n-p}}\bigg) \\ &=& \frac{1}{(\varphi_1)'_l(1)} \cdot \varphi_2\bigg(\bigg(\frac{C_p(\widetilde{\Omega})}{C_p(\Omega)}\bigg)^{\frac{1}{n-p}}\bigg),
\end{eqnarray*} 
where $(\varphi_1)'_l(1)$ is assumed to exist and to be nonzero. Together with inequality (\ref{C-O-B-M-inequality-2-111}), one gets
\begin{equation*}  
   (\varphi_1)'_l(1)\cdot \lim_{\varepsilon\rightarrow 0^+} \frac{1- \left( \frac{C_p(\Omega)}  {C_p(\Omega+_{\varphi, \varepsilon} \widetilde{\Omega})}\right)^{\frac{1}{n-p}}}  {\varepsilon}\geq   \varphi_2\bigg(\bigg(\frac{C_p(\widetilde{\Omega})}{C_p(\Omega)}\bigg)^{\frac{1}{n-p}}\bigg).
\end{equation*}  
This together with Theorem \ref{c-interpretation-111} further imply the $p$-capacitary Orlicz-Minkowski inequality:
\begin{eqnarray*}  
 	(n-p)\cdot C_{p, \varphi_2}(\Omega, \widetilde{\Omega})   &=& (\varphi_1)'_l(1)\cdot \lim_{\varepsilon\rightarrow 0^+} C_p(\Omega +_{\varphi, \varepsilon}\widetilde{\Omega}) \cdot \lim_{\varepsilon\rightarrow 0^+} \frac{1-\frac{C_p(\Omega)}  {C_p(\Omega+_{\varphi, \varepsilon} \widetilde{\Omega})}}  {\varepsilon} \\ & =&(\varphi_1)'_l(1)\cdot  (n-p)\cdot   C_p(\Omega) \cdot \lim_{\varepsilon\rightarrow 0^+} \frac{1- \left( \frac{C_p(\Omega)}  {C_p(\Omega+_{\varphi, \varepsilon}\widetilde{\Omega})}\right)^{\frac{1}{n-p}}}  {\varepsilon}\\ &\geq&  (n-p)\cdot C_p(\Omega) \cdot   \varphi_2\bigg(\bigg(\frac{C_p(\widetilde{\Omega})}{C_p(\Omega)}\bigg)^{\frac{1}{n-p}}\bigg). 
\end{eqnarray*}

On the other hand, assume that  the $p$-capacitary  Orlicz-Minkowski inequality in Theorem \ref{capacitary-Minkowski-2-2-1-1} holds. In particular, for $\varphi_1, \varphi_2\in \Phi_1$ and for $\Omega, \widetilde{\Omega}\in \cC_0$,  
\begin{eqnarray*}   
	\frac{C_{p, \varphi_1}(+_{\varphi}(\Omega, \widetilde{\Omega}), \Omega)}{C_p(+_{\varphi}(\Omega, \widetilde{\Omega}))}  &\geq& \varphi_1\bigg(\bigg(\frac{C_p(\Omega)}{C_p(+_{\varphi}(\Omega, \widetilde{\Omega}))}\bigg)^{\frac{1}{n-p}}\bigg), \\ \frac{C_{p, \varphi_2}(+_{\varphi}(\Omega, \widetilde{\Omega}), \widetilde{\Omega})}{C_p(+_{\varphi}(\Omega, \widetilde{\Omega}))}  &\geq&  \varphi_2\bigg(\bigg(\frac{C_p(\widetilde{\Omega})}{C_p(+_{\varphi}(\Omega, \widetilde{\Omega}))}\bigg)^{\frac{1}{n-p}}\bigg),
\end{eqnarray*}  
where  $\varphi=\alpha_1\varphi_1+\alpha_2\varphi_2$ with $\alpha_1, \alpha_2>0$ and $\varphi_1, \varphi_2\in \Phi_1$,  and $+_{\varphi}(\Omega, \widetilde{\Omega})$ is the convex domain whose support function $h_{+_{\varphi}(\Omega, \widetilde{\Omega})}$ is given by    
\begin{eqnarray*}  
	1 = \alpha_1 \varphi_1\bigg(\frac{h_{\Omega}(u)}{h_{+_{\varphi}(\Omega, \widetilde{\Omega})}(u)}\bigg)+\alpha_2 \varphi_2\bigg(\frac{h_{\widetilde{\Omega}}(u)}{h_{+_{\varphi}(\Omega, \widetilde{\Omega})}(u)}\bigg), \ \ \ \mathrm{for} \ \ u\in \sphere.
\end{eqnarray*} 
This together with (\ref{def:mixed:orlicz:capacity}) lead to inequality (\ref{Orlicz-Minkowski-inequality-linear--1}) with $m=2$: 
\begin{eqnarray*}  
	1\!\! &=&\!\! \frac{p-1}{n-p}\cdot \int_{\sphere}\bigg[\alpha_1 \varphi_1\bigg(\frac{h_{\Omega}(u)}{h_{+_{\varphi}(\Omega, \widetilde{\Omega})}(u)}\bigg)+\alpha_2 \varphi_2\bigg(\frac{h_{\widetilde{\Omega}}(u)}{h_{+_{\varphi}(\Omega, \widetilde{\Omega})}(u)}\bigg)\bigg]\cdot \frac{ h_{ +_{\varphi}(\Omega, \widetilde{\Omega})}(u)}{C_p(+_{\varphi}(\Omega, \widetilde{\Omega}))} \cdot \,d\mu_p(+_{\varphi}(\Omega, \widetilde{\Omega}), u) \\ \!\! &=&\!\! \alpha_1\cdot \frac{C_{p, \varphi_1}(+_{\varphi}(\Omega, \widetilde{\Omega}), \Omega)}{C_p(+_{\varphi}(\Omega, \widetilde{\Omega}))}+\alpha_2 \cdot \frac{C_{p, \varphi_2}(+_{\varphi}(\Omega, \widetilde{\Omega}), \widetilde{\Omega})}{C_p(+_{\varphi}(\Omega, \widetilde{\Omega}))} \\ \!\! &\geq&\!\!   \alpha_1\cdot \varphi_1\bigg(\bigg(\frac{C_p(\Omega)}{C_p(+_{\varphi}(\Omega, \widetilde{\Omega}))}\bigg)^{\frac{1}{n-p}}\bigg)+\alpha_2\cdot  \varphi_2\bigg(\bigg(\frac{C_p(\widetilde{\Omega})}{C_p(+_{\varphi}(\Omega, \widetilde{\Omega}))}\bigg)^{\frac{1}{n-p}}\bigg).
\end{eqnarray*}

The $M$-addition of convex domains are closely related to the Orlicz addition. For an arbitrary subset $M\subset \R ^m$,  the $M$-addition of $\Omega_1, \cdots, \Omega_m\in \cC_0$, denoted by $\oplus_M(\Omega_1,\cdots, \Omega_n)$, is defined by (see e.g. \cite{GHW2013, GHW2014, Pro99, Pro97})
$$
   \oplus_M(\Omega_1, \cdots, \Omega_m)=\bigg\{\sum_{j=1}^{m}a_jx^j: \ \ x^j\in \Omega_j \ \ \mathrm{and}\ \  (a_1, \cdots, a_m)\in M \bigg\}.
$$
It is equivalent to the following more convenient formula:
\begin{equation}\label{m addition}
   \oplus_M(\Omega_1, \cdots, \Omega_m)=\cup\big\{a_1 \Omega_1+\cdots +a_m\Omega_m: \ \ (a_1, a_2,\cdots, a_m)\in M \big\},
\end{equation} 
where $a_1 \Omega_1+\cdots +a_m\Omega_m$ is the Minkowski addition of $a_j\Omega_j=\{a_jx^j: x^j\in \Omega_j\}$ for $j=1, 2, \cdots, m$.  Note that if $M$ is compact, then $\oplus_M (\Omega_1, \cdots, \Omega_m)$ is again a convex domain. In general, the $M$-addition is different from the Orlicz addition. However,  when $M$ is a $1$-unconditional convex body in $\R^m$ that contains $\{e_1,\cdots, e_m\}$ in its boundary, then the $M$-addition coincides with the Orlicz $L_{\varphi}$ addition for some $\varphi\in\Phi_m$. More properties and historical remarks for the $M$-addition, such as convexity, $GL(n)$ covariance,  homogeneity and monotonicity, can be founded in \cite{GHW2013, GHW2014, Pro99, Pro97}.

\bl\label{equality condition} 
If $M\subset \R ^m$ is compact and $\Omega_1,\cdots, \Omega_m \in \cC_0$, then for any $a=(a_1,\cdots,a_m)\in M$, 
\begin{equation} \label{inequality for m addition}
   C_p\big(\!\oplus_M(\Omega_1,\cdots, \Omega_m)\big)^{\frac{1}{n-p}}\geq\sum_{i=1}^{m}\bigg[|a_i|\cdot C_p(\Omega_i)^{\frac{1}{n-p}}\bigg]. 
\end{equation} 
If equality holds in (\ref{inequality for m addition}) for some $a\in M$ with $a_j\neq 0$ for all  $j=1, 2, \cdots, m$,  then $\Omega_i$ is homothetic to $\Omega_j$ for all $1\leq i<j\leq m$.
\el
\begin{proof} 
Recall that the $p$-capacity is invariant under affine isometries and has homogeneous degree $n-p$ (see \cite{EvansGariepy}). Then for all $a\in \R$ and for all $\Omega\in \cC_0$, one has 
$$
   C_p(a\Omega)=|a|^{n-p} C_p(\Omega).
$$ 
Note that $n-p>0$. It follows from (\ref{Brunn-Minkowski-classical-p-capacity}), (\ref{m addition}) and the monotonicity of the $p$-capacity that, for all  $a=(a_1,\cdots,a_m)\in M$,  
$$ 
   C_p\big(\!\oplus_M(\Omega_1,\cdots, \Omega_m)\big)^{\frac{1}{n-p}} \geq C_p \big(a_1\Omega_1+\cdots +a_m\Omega_m \big)^{\frac{1}{n-p}}\geq \sum_{i=1}^{m}\bigg[|a_i|\cdot C_p(\Omega_i)^{\frac{1}{n-p}}\bigg].
$$
Assume that equality holds in (\ref{inequality for m addition}) for some $a\in M$ with $a_j\neq 0$ for all $j=1, 2, \cdots, m$. Then equality in (\ref{Brunn-Minkowski-classical-p-capacity}) must hold and hence $\Omega_i$ is homothetic to $\Omega_j$ for all $1\leq i<j\leq m$.  \end{proof}
   
Let $e_j^{\perp}=\{x\in \R^m: \langle x, e_j\rangle=0\}$ for all $j=1, 2, \cdots, m$.  For a nonzero vector $x\in \R^m$ and a convex set $E\subset \R^m$, define the support set of $E$ with outer normal vector $x$ to be the set 
$$
   F(E, x)=\left\{y\in \R^m: \langle x, y\rangle =\sup_{z\in E} \langle x, z\rangle \right\} \cap E.
$$ 
\bt 
Let  $M\subset \R ^m$ be a compact subset and $\Omega_1,\cdots, \Omega_m\in \cC_0$. Then
\begin{equation}\label{inequality 1 for m addition}
   C_p\big(\!\oplus_M(\Omega_1,\cdots, \Omega_m)\big)^{\frac{1}{n-p}} \geq h_{conv( M)} \big(C_p(\Omega_1)^{\frac{1}{n-p}}, \cdots, C_p(\Omega_m)^{\frac{1}{n-p}}\big).
\end{equation}
If $M\cap F(conv(M), x) \not\subset \cup_{j=1}^m e_j^{\perp}$ for all $x=(x_1, \cdots, x_m)$ with all $x_i>0$ and equality holds in (\ref{inequality 1 for m addition}), then  $\Omega_i$ is homothetic to $\Omega_j$ for all $1\leq i<j\leq m$. 
\et
\begin{proof} 
It is easily checked that  $h_{conv(M)}(x)=\max_{y\in M} \langle x, y\rangle$ for all $x\in \R^m$. Following (\ref{inequality for m addition}), one has, as all $C_p(\Omega_i)>0$, 
\begin{eqnarray*} 
   C_p\big(\!\oplus_M(\Omega_1,\cdots, \Omega_m)\big)^{\frac{1}{n-p}} &\geq & \max_{(a_1, \cdots, a_m)\in M}\big \langle \big(|a_1|, \cdots, |a_m|\big), \big(C_p(\Omega_1)^{\frac{1}{n-p}}, \cdots, C_p(\Omega_m)^{\frac{1}{n-p}}\big) \big \rangle \\ &\geq & \max_{(a_1, \cdots, a_m)\in M}\big \langle \big(a_1, \cdots, a_m\big), \big(C_p(\Omega_1)^{\frac{1}{n-p}}, \cdots, C_p(\Omega_m)^{\frac{1}{n-p}}\big) \big \rangle \\ &=& h_{conv( M)} \big(C_p(\Omega_1)^{\frac{1}{n-p}}, \cdots, C_p(\Omega_m)^{\frac{1}{n-p}}\big).  
\end{eqnarray*} 
Now let us characterize the conditions for equality. Let $$x_0=\big(C_p(\Omega_1)^{\frac{1}{n-p}}, \cdots, C_p(\Omega_m)^{\frac{1}{n-p}}\big).$$  Assume that equality holds in (\ref{inequality 1 for m addition}). There exists a vector $a_0\in M\cap F(conv(M), x_0)$ such that  
\begin{eqnarray*} 
   C_p\big(\!\oplus_M(\Omega_1,\cdots, \Omega_m)\big)^{\frac{1}{n-p}} &= & \max_{(a_1, \cdots, a_m)\in M}\big \langle \big(|a_1|, \cdots, |a_m|\big), \big(C_p(\Omega_1)^{\frac{1}{n-p}}, \cdots, C_p(\Omega_m)^{\frac{1}{n-p}}\big) \big \rangle \\ &= & \max_{(a_1, \cdots, a_m)\in M}\big \langle \big(a_1, \cdots, a_m\big), \big(C_p(\Omega_1)^{\frac{1}{n-p}}, \cdots, C_p(\Omega_m)^{\frac{1}{n-p}}\big) \big \rangle \\ &=& h_{conv( M)} (x_0)=\langle a_0, x_0\rangle.  
\end{eqnarray*} 
Note that $M\cap F(conv(M), x_0) \not\subset \cup_{j=1}^m e_j^{\perp}$ and then  all coordinates of $a_0$ must be strictly positive. As all coordinates of $x_0$ are strictly positive, it follows from the conditions of equality for (\ref{inequality for m addition}) that  $\Omega_i$ is homothetic to $\Omega_j$ for all $1\leq i<j\leq m$.  
\end{proof}

\section{The $p$-capacitary Orlicz-Minkowski problems} 

For $\varepsilon\in \R$ close to $0$, $\varphi\in \mathcal{D}\cup\mathcal{I}$, and  a continuous function $\psi: (0, \infty)\rightarrow (0, \infty)$, consider the function  $\widehat{f}_{\varepsilon}: \sphere \rightarrow (0, \infty)$  defined by \begin{eqnarray}\label{new:Orlicz:addition} 
   \widehat{f}_{\varepsilon}(u)= {\varphi}^{-1}\Big(\varphi\big(h_{\Omega}(u)\big)+\varepsilon \psi\big(h_{\Omega_1}(u)\big)\Big) 
\end{eqnarray}   
for $\Omega, \Omega_1\in \cC_0$.  Clearly, there exists $\varepsilon_0>0$ such that  $\widehat{f}_{\varepsilon}\in C^+(\sphere)$ for all $\varepsilon\in (-\varepsilon_0, \varepsilon_0)$.  Similar to Lemmas 6 and 7 in \cite{HLYZ}, if the derivative of $\varphi$ (denoted by $\varphi'$) exists, and is strictly positive and continuous on $(0, \infty)$, then the following limit is uniform on $\sphere$:  \begin{equation}\label{conv-hat-sum}  
   \lim_{\varepsilon\rightarrow 0} \frac{\widehat{f}_\varepsilon(u)-h_{\Omega}(u)}{\varepsilon}= \frac{\psi\left(h_{\Omega_1}(u)\right)}{\varphi'\left(h_{\Omega}(u)\right)} 
\end{equation} 
and hence 
$\widehat{f}_\varepsilon$ converges uniformly to $h_{\Omega}$ on $\sphere$.  For completeness, a brief proof of (\ref{conv-hat-sum}) is presented here. By the chain rule and $(\varphi^{-1})'(a)=1/\varphi'(\varphi^{-1}(a))$, one has,  on $\sphere\times (-\varepsilon_0, \varepsilon_0)$,  $$\frac{\partial \widehat{f}_\varepsilon}{\partial  \varepsilon} \big |_{(u, \varepsilon)}=\frac{\psi\left(h_{\Omega_1}(u)\right)}{\varphi' (
 \widehat{f}_\varepsilon)}\ \ \ \mathrm{and}\ \ \  \frac{\partial \widehat{f}_\varepsilon}{\partial  \varepsilon} \big |_{(u, 0)}= \frac{\psi\left(h_{\Omega_1}(u)\right)}{\varphi'\left(h_{\Omega}(u)\right)}.$$  It can be checked that $\partial \widehat{f}_\varepsilon /\partial  \varepsilon$ is continuous on $\sphere\times (-\varepsilon_0, \varepsilon_0)$. Hence, $\partial \widehat{f}_\varepsilon /\partial  \varepsilon$ is uniformly continuous,  has strictly positive lower bounds, and has finite upper bounds on the compact set $\sphere\times  [-\varepsilon_0/2, \varepsilon_0/2]$. For every $\varepsilon\in [-\varepsilon_0/2, \varepsilon_0/2]$, it follows from  the mean value theorem that 
$$
   \bigg|  \frac{\widehat{f}_\varepsilon(u)-h_{\Omega}(u)}{\varepsilon}-\frac{\partial \widehat{f}_\varepsilon}{\partial  \varepsilon} \big |_{(u, 0)} \bigg|=\bigg|  \frac{\partial \widehat{f}_\varepsilon}{\partial  \varepsilon} \big |_{(u, \varepsilon')}-\frac{\partial \widehat{f}_\varepsilon}{\partial  \varepsilon} \big |_{(u, 0)} \bigg|
$$ 
for some $\varepsilon'\in [-\varepsilon, \varepsilon]\subset  [-\varepsilon_0/2, \varepsilon_0/2]$. The uniform convergence of  (\ref{conv-hat-sum}) then follows from the uniform continuity of  $\partial \widehat{f}_\varepsilon /\partial  \varepsilon$ on $\sphere\times  [-\varepsilon_0/2, \varepsilon_0/2]$.  
 
 An argument similar to the proof of Theorem \ref{c-interpretation-111} yields the following result regarding the asymptotic behavior of $C_p(\Omega_{\widehat{f}_\varepsilon})$, the $p$-capacity of the Aleksandrov domain associated to $\widehat{f}_{\varepsilon}$.  Similar arguments for the volume can be found in \cite{GHW2014, HLYZ}. 

\bp\label{p-5-anther-mixed} 
Let $\varphi\in \mathcal{I}\cup \mathcal{D}$ be such that $\varphi'$ exists, and is nowhere zero and continuous on $(0, \infty)$. For $\Omega, \Omega_1\in \cC_0$, one has
$$
   \frac{1}{n-p}\cdot \lim_{\varepsilon \rightarrow 0}\frac{C_p(\Omega_{\widehat{f}_\varepsilon})-C_p(\Omega)}{\varepsilon}=\frac{p-1}{n-p} \cdot \int_{S^{n-1}}
   \frac{\psi\left(h_{\Omega_1}(u)\right)}{\varphi'\left(h_{\Omega}(u)\right)}\,d\mu_p(\Omega, u).
$$  
\ep

\noindent {\bf Remark.} Assume that $\phi=1/\varphi'$ and then \begin{equation}\label{relation of two functions}
   \varphi(t)=\int_0^t \frac{1}{\phi(s)}\,ds.
\end{equation} 
Proposition \ref{p-5-anther-mixed} can be rewritten as 
$$
   \frac{1}{n-p}\cdot \lim_{\varepsilon \rightarrow 0}\frac{C_p(\Omega_{\widehat{f}_\varepsilon})-C_p(\Omega)}{\varepsilon}=\frac{p-1}{n-p} \cdot \int_{S^{n-1}}
   \psi\left(h_{\Omega_1}(u)\right)\cdot  \phi\left(h_{\Omega}(u)\right) \,d\mu_p(\Omega, u).
$$ 
While if $\varphi\in \mathcal{D}$ and $\phi=-1/\varphi'$, then 
\begin{equation*}
   \varphi(t)=\int_{t}^{\infty} \frac{1}{\phi(s)}\,ds. 
\end{equation*}  
Proposition \ref{p-5-anther-mixed} can be rewritten as 
$$
   \frac{1}{n-p}\cdot \lim_{\varepsilon \rightarrow 0}\frac{C_p(\Omega)-C_p(\Omega_{\widehat{f}_\varepsilon})}{\varepsilon}=\frac{p-1}{n-p} \cdot \int_{S^{n-1}}
   \psi\left(h_{\Omega_1}(u)\right)\cdot  \phi\left(h_{\Omega}(u)\right) \,d\mu_p(\Omega, u).
$$

It is not clear whether there are $p$-capacitary Orlicz-Brunn-Minkowski and Orlicz-Minkowski inequalities involving the addition defined by (\ref{new:Orlicz:addition}). It is worth to mention that Proposition \ref{p-5-anther-mixed} provides a geometric meaning of the measure 
$$
   \mu_{p, \phi}(\Omega, \cdot)=\phi(h_{\Omega})\,d\mu_p(\Omega, \cdot).
$$
The measure  $\mu_{p, \phi}(\Omega, \cdot)$ will be called the Orlicz $L_{\phi}$ $p$-capacitary measure of $\Omega$. When $\Omega\in \cC_0$ and $\phi: (0, \infty)\rightarrow (0, \infty)$ is a continuous function, $\phi(h_{\Omega})$ has strictly positive lower bound and finite upper bound on $\sphere$. This further implies that the support of the measure $\mu_{p, \phi}(\Omega, \cdot)$ is not contained in any closed hemisphere and  by (\ref{non-concentration--1}) 
\begin{equation*}  
   \int_{\sphere} \langle \theta, u\rangle_+ \,d\mu_{p, \phi}(\Omega, u)> 0 \ \ \ \ \ \mathrm{for\ each} \ \theta\in \sphere.
\end{equation*}  
Although the $p$-capacity is translation invariant, one {\it cannot} expect to have the centroid of $\mu_{p, \phi}(\Omega, \cdot)$ at the origin even if $\phi=t^{1-q}$ for all $q\neq 1$ (as the $L_q$ sum is not linear). When $\phi=t^{1-q}$ (and hence $\varphi=t^q/q$) for $q\neq 0$,  one gets the $L_q$ $p$-capacitary measure of $\Omega$ which will be denoted by $\mu_{p, q}(\Omega, \cdot)$.

It is interesting and important to study the following $p$-capacitary Orlicz-Minkowski problem: {\em given a fixed continuous function $\phi: (0, \infty)\rightarrow (0, \infty)$ and a finite Borel measure $\mu$ on $\sphere$, does there exist a convex domain $\Omega$ whose closure $\overline{\Omega}$ contains the origin such that 
$$
    \frac{\mu}{\phi(h_{\Omega})}=\tau \cdot \mu_p(\Omega,\cdot) \ \ (\mathrm{or} \  \  \mu= \tau\cdot \mu_{p, \phi}(\Omega, \cdot) \ \ \mathrm{if\ } \  \Omega\in \cC_0)  
$$ 
for some positive number $\tau$?}  When  $\phi=t^{q-1}$ for $q\neq 0$, we are interested in the following normalized  $p$-capacitary $L_q$ Minkowski problem:  {\em given a finite Borel measure $\mu$ on $\sphere$, does there exist a convex domain $\Omega$ whose closure $\overline{\Omega}$ contains the origin such that 
$$
   \mu \cdot h^{q-1}_{\Omega}=c_q \cdot \frac{\mu_p(\Omega,\cdot)}{C_p(\Omega)} \ \ \ (\mathrm{or}\ C_p(\Omega)\cdot  \mu  = c_q \cdot   \mu_{p, q} (\Omega,\cdot) \ \ \mathrm{if}\ \  \Omega\in \cC_0). 
$$  
for some positive number $c_q$?}  The $p$-capacitary $L_1$ Minkowski problem (i.e., $q=1$) has been studied in \cite{AGHLV, CNSXYZ,  Jerison, Jerison-1996}.   

In this section, we will provide a solution for the $p$-capacitary Orlicz-Minkowski problem as well as the normalized $p$-capacitary $L_q$ Minkowski problem for $q>1$. 
Throughout the rest of this section, unless otherwise stated, let $(\phi, \varphi)$ be the pair of functions such that 
\begin{itemize}   [topsep=1pt,itemsep=-0.5ex,partopsep=1ex,parsep=1ex, leftmargin=.6in]
\item [(A1):]  the function $\phi:(0,\infty)\rightarrow(0,\infty)$ is decreasing and continuous with  $\lim_{t\rightarrow 0^+} \phi(t)=\infty$,  
\item  [(A2):] the function $\varphi$ given by (\ref{relation of two functions}) satisfies $\varphi(t)<\infty$ for all $t>0$ and $\lim_{t\rightarrow \infty} \varphi(t)=\infty$.
\end{itemize}
Obviously $\varphi$ is strictly increasing such that  $
\varphi(0)= \lim_{t\rightarrow 0^+}\varphi(t)=0$ and $\lim_{t\rightarrow 0^+}\varphi'(t)=0.$ The inverse of $\varphi$, denoted by $\varphi^{-1}$, exists and is also continuously differentiable on $(0,\infty)$.  Moreover $(\varphi^{-1})'(t)=\phi(\varphi^{-1}(t))$ for all $t\in (0, \infty)$,   
$\varphi^{-1}(0)= \lim_{t\rightarrow 0^+}\varphi^{-1}(t)=0$ and $\lim_{t\rightarrow \infty}\varphi^{-1}(t)=\infty.$ It is easily checked that for all $a,b>0$, there are constants $t_0, M_1, M_2>0$, such that, for all $t\in (0, t_0)$ (see e.g., \cite[Lemma 4.1(iii)]{huang}), 
$$
   M_1\leq (\varphi^{-1})'(\varphi(a)-b\varphi(t))\leq M_2. 
$$    
For convenience, we use $\|f\|_{\varphi, \mu}$ to denote the ``Orlicz norm" of $f\in C(\sphere)$, where $C(\sphere)$ denotes the set of all continuous functions on $\sphere$: 
$$
   \|f\|_{\varphi, \mu}=\inf\Big\{ \lambda>0: \int _{\sphere} \varphi\Big(\frac{f}{\lambda}\Big)\,d\mu \leq \varphi(1)\cdot \int _{\sphere} \,d\mu\Big\},
$$ 
where $(\phi, \varphi)$ satisfies conditions (A1) and (A2). Clearly, $\|f\|_{\varphi, \mu} \leq \|g\|_{\varphi, \mu}$  if $f\leq g$ and $\|af\|_{\varphi, \mu}=a\|f\|_{\varphi, \mu}$ for all $a\geq 0$.

\subsection{The $p$-capacitary Orlicz-Minkowski problem of discrete measures} \label{section:polytope-minkowski-solution} 

In this subsection,  we provide a solution for the $p$-capacitary Orlicz-Minkowski problem for discrete measure $\mu$ under the very limit condition: {\em  the support of $\mu$ is not contained in any closed hemisphere and either $\mu(\{\xi\})=0$ or $\mu(\{-\xi\})=0$ for all $\xi\in S^{n-1}$.} For simplicity, we always use $P'$ to denote the interior of a polytope $P$. 

\bt\label{discrete solution----1}  
Suppose that $(\phi, \varphi)$ satisfies conditions (A1) and (A2).  Let $\mu=\sum_{i=1}^{m}\lambda_i\delta_{u_i}$ be such that either $\mu(\{\xi\})=0$ or $\mu(\{-\xi\})=0$ for all $\xi\in S^{n-1}$,  where $\lambda_1,\lambda_2,\cdots, \lambda_m>0$ are given constants and $\{u_1,u_2,\cdots, u_m\} \subset S^{n-1}$ are not contained in any closed hemisphere. There exists a polytope $P$ with the origin in its interior and a constant $\tau>0$ such that for $1<p<n$, 
$$
   \mu=\tau\cdot \phi(h_P) \cdot \mu_p(P', \cdot)=\tau \cdot \mu_{p, \phi}(P', \cdot)\ \ \ \mathrm{or}\ \ \ \frac{\mu}{\phi(h_P)}=\tau\cdot \mu_p(P', \cdot).
$$ 
Moreover, $\|h_P\|_{\varphi, \mu}=1$ and the constant $\tau$ can be calculated by 
$$
   \tau=\frac{\varphi(1) \cdot \int_{\sphere}\,d\mu(u)}{\int_{\sphere} \varphi(h_{P}(u))\,d\mu_{p, \phi}(P', u)}= \bigg(\frac{p-1}{n-p} \bigg)\cdot \frac{1}{C_p(P')} \cdot\int_{\sphere} \frac{h_P(u)}{\phi(h_P(u))}\,d\mu(u).
$$
\et

In order to prove Theorem \ref{discrete solution----1}, we need the following lemma. For $u\in \sphere$ and $t\geq 0$, let $H_{u,t}^{-}=\big\{y\in \R^n:  \ \langle u,  y\rangle \leq t\big\}.$ Let $m>n$ be an integer and  
$$
  \R_{*}^{m}=\{(x_1, \cdots, x_m): \  x_i\geq 0\}\ \ \mathrm{and} \ \  \R_{+}^{m}=\{(x_1, \cdots, x_m): \ x_i>0\}.
$$
Define the polytope $P(x)$ for $x\in\R_{+}^{m}$ by  $P(x)=\cap_{i=1}^{m}H_{u_i,x_i}^{-}$ with $\{u_1, \cdots, u_m\}\subset \sphere$.    The following lemma states the differentiability of $C_p(P(x))$ (see \cite[Lemma 3.2]{HugLYZ} for the volumetric analogue).

\bl\label{derivative of capacity} 
Let $p\in (1, n)$. Suppose that the vectors $u_1,u_2,\cdots,u_m\in S^{n-1}$ are not contained in any closed hemisphere. {Then $C_p(P'(x))$} is differentiable and  for $1\leq i\leq m$,  
$$
   \frac{\partial}{\partial x_i} C_p(P'(x))=(p-1) \cdot \mu_p(P'(x), \{u_i\}).
$$ 
\el

\begin{proof} It follows from \cite[Theorem 5.2]{CNSXYZ} that
for $f\in C^{+}(S^{n-1})$ and $g\in C(\sphere)$, 
$$
   \frac{\,dC_p(\Omega_{f+\varepsilon g})}{\,d\varepsilon}\big|_{\varepsilon=0}=(p-1)\int_{S^{n-1}}g(u)\, d\mu_p(\Omega_{f}, u),
$$
where $\Omega_{f+\varepsilon g}$ and $\Omega_f$ are the Aleksandrov domains associated to $f+\varepsilon g$ and $f$, respectively.  Let $g_i\in C(\sphere)$ for all $i=1, 2, \cdots, m$ be  such that $g_i(u_i)=1$ and $g_i(u_j)=0$ if $i\neq j$. Let $f=h_{P(x)}$ and then  $h_{P(x)}(u_j)=x_j$ for $j=1, 2, \cdots, m$. Thus
\begin{eqnarray*}  
	\frac{\partial}{\partial x_i} C_p(P'(x))
    &=& \lim_{t\rightarrow 0}\frac{C_p(P'(x_1, \cdots, x_{i-1}, x_i+t, x_{i+1}, \cdots, x_m))-C_p(P'(x))}{t}\\
    &=& \lim_{t\rightarrow 0}\frac{C_p(\Omega_{f+tg_i})-C_p(\Omega_f)}{t}\\
    &=& (p-1)\int_{S^{n-1}}g_i(u)\ d\mu_p(P'(x),u)\\
    &=& (p-1)\int_{\{u_1,\cdots,u_m\}}g_i(u)\ d\mu_p(P'(x),u)\\
    &=& (p-1)\cdot \mu_p(P'(x), \{u_i\})
\end{eqnarray*}
and the desired result is obtained. 
\end{proof}

\noindent {\em Proof of Theorem \ref{discrete solution----1}.} Our proof is based on the techniques in \cite{huang, HugLYZ}.   Let $\mu=\sum_{i=1}^{m}\lambda_i\delta_{u_i}$ be the given finite discrete Borel measure on $\sphere$ where $\lambda_1,\lambda_2,\cdots, \lambda_m>0$ are given constants and $\{u_1,u_2,\cdots, u_m\} \subset S^{n-1}$ are not contained in any closed hemisphere. Let \begin{equation}\label{construct:polytopes}
P(x)=\bigcap_{i=1}^{m}H_{u_i,x_i}^{-}\end{equation} for $x=(x_1, \cdots, x_m)\in \R^m_*$. Consider the optimization problem: $\sup_{x\in M^*} C_p(P(x))$ where $M^*$ is the compact surface in $\R^m$:  
$$
   M^*=\Big\{x\in \R_{*}^{m}: \ \ \sum_{i=1}^{m}\lambda_i\varphi(x_i)= \varphi(1) \cdot \sum_{i=1}^m \lambda_i  \Big\}.
$$ 
Note that $P(x)$ for $x\in M^*$ defines a compact convex set because $\{u_1,u_2,\cdots, u_m\} \subset S^{n-1}$ are not contained in any closed hemisphere. In fact, as $M^*$ is compact, then $x\in M^*$ must have $x_i<\infty$ for all $i=1, \cdots, m$. On the other hand, by (\ref{construct:polytopes}) and $\{u_1, \cdots, u_m\}$ positively spans $\Rn$ (see \cite[p.411]{SchneiderBook}), $P(x)$ is circumscribed by the hyperplanes $H_{u_i,x_i}=\{z\in \Rn: \langle z, u_i\rangle =x_i\}$ and hence $P(x)$ for any $x\in M^*$ is bounded. Of course,  $o\in P(x)$ and $x_i\geq h_{P(x)}(u_i)$ for all $1\leq i\leq m$ and all $x\in M^*$,  with $x_i= h_{P(x)}(u_i)$  if  $S(P(x), \{u_i\})>0$ (and hence $\mu_p(P'(x), \{u_i\})>0$, due to Lemma \ref{relation between capacity measure and surface area measure}). Moreover, $C_p(P(x))$ is continuous about $x\in  \R_{*}^{m}$ due to the continuity of the $p$-capacity and the continuity of $P(x)$ about $x\in \R_{*}^{m}$ with respect to the Hausdorff metric (see  \cite[page 57]{SchneiderBook}). As $M^*$ is compact, there exists $z\in M^*$ such that  
$$
   C_p(P(z))=\sup_{x\in M^*}  C_p(P(x)).
$$

By \cite[Proposition 13.2]{AGHLV}, an argument similar to the one in \cite[Section 13.1]{AGHLV}, where $\hat{q}_i$, $E_1$ and $k(t)$ in \cite{AGHLV} are replaced by $z_i$,  $P(z)$ and $\sum_{i=1}^{m}\lambda_i \varphi\Big(\Big(\frac{C_p(P(z))}{C_p(E_1+tE_2)}\Big)^{\frac{1}{n-p}}\cdot (z_i+at)\Big),$  yields that $P(z)$  has nonempty interior. 
Of course,  $o\in P(z)$ and $z_i\geq h_{P(z)}(u_i)$ for all $1\leq i\leq m$,  with $z_i= h_{P(z)}(u_i)$  if  $S(P(z), \{u_i\})>0$ (and hence $\mu_p(P(z), \{u_i\})>0$, due to Lemma \ref{relation between capacity measure and surface area measure}).

It remains to show that $o\in P'(z)$. To this end, assume that $o\in \partial P'(z)$. For simplicity, let $h_i=h_{P'(z)}(u_i)$ for $1\leq i\leq m$. Without loss of generality, let $h_1=\cdots=h_k=0$ and $h_{k+1},\cdots,h_m >0$ for some $1\leq k<m$. The fact that $k<m$ follows from $B^n_2\subset P(x_0)$  and   \begin{eqnarray}\label{uniform-lower-bound} 
   C_p(P'(z))\geq C_p(P'(x_0))\geq C_p(B^n_2)>0, 
\end{eqnarray}  
where $B^n_2$ is the unit open ball in $\Rn$ and $x_0=(1, 1, \cdots, 1)\in M^*$. 

In order to get a contradiction with the maximality of $C_p(P'(z))$, we need to construct $z^t\in M^*$ such that $C_p(P'(z^t))>C_p(P'(z))$. For $t>0$ small enough, let $z^t=(z_1^t, \cdots, z_m^t)\in M^*$ be given by \begin{equation*}
   z^{t}_i=
    \begin{cases}
      \varphi^{-1}(\varphi(z_i)+\varphi(t)), & 1\leq i\leq k \\
      \varphi^{-1}(\varphi(z_i)-\lambda\varphi(t)), & k+1\leq i\leq m,
    \end{cases}
\end{equation*}
where  $\lambda=\frac{\lambda_1+\cdots +\lambda_k}{\lambda_{k+1}+\cdots+\lambda_m}.$ For $ k+1\leq i\leq m$, let $h^{t}_{i}=\varphi^{-1}(\varphi(h_i)-\lambda\varphi(t))$ and then \begin{eqnarray*} 
   \lim_{t \rightarrow 0^+}\frac{h_{i}^{t}-h_i}{t}&=&\lim_{t\rightarrow 0^+}(\varphi^{-1})'(\varphi(h_i)-\lambda\varphi(t))(-\lambda\varphi'(t))=0, \\   h^{t}_{i} &\leq&  \varphi^{-1}(\varphi(z_i)-\lambda\varphi(t))=z_i^t, 
\end{eqnarray*}  
where the inequality is due to $h_i=h_{P(z)}(u_i)\leq z_i$  and the fact that  both $\varphi$ and $\varphi^{-1}$ are strictly increasing. Similarly $t\leq z^{t}_i$ for all $1\leq i\leq k$  and hence 
$$
   P^t=\Big(\bigcap_{i=1}^{k}H^{-}_{u_i, t}\Big) \bigcap \Big(\bigcap^{m}_{i=k+1}H^{-}_{u_i,h^{t}_{i}}\Big)\subset P(z^t).
$$ 
Note that $(P^0)'=P'(z)$ and  $o\in (P^t)'$ if $t>0$ is small enough. Moreover,
\begin{eqnarray*} 
   C_p((P^t)') &=&\frac{p-1}{n-p}\cdot\bigg( t\sum_{i=1}^{k}\mu_p((P^t)', \{u_i\}) +\sum^{m}_{i=k+1}h_{i}^{t}\mu_p((P^t)', \{u_i\})\! \bigg), \\   C_p((P^t)',P'(z)) &=& \frac{p-1}{n-p}\cdot \sum_{i=k+1}^{m}h_i\mu_p((P^t)', \{u_i\}). 
\end{eqnarray*} 
As $t\rightarrow 0^+$, one has $P^t\rightarrow P(z)$ in the Hausdorff metric and $\mu_p((P^t)', \cdot)$ converges to $\mu_p(P'(z), \cdot)$ weakly  (see \cite[Lemma 4.1]{CNSXYZ}). Moreover, 
\begin{eqnarray*}
   \lim_{t\rightarrow 0^+}\frac{C_p((P^t)')-C_p((P^t)', P'(z))}{t} &=&
   \frac{p-1}{n-p}\lim_{t \rightarrow 0^+}\bigg(\!\sum_{i=1}^{k} \mu_p((P^t)', \{u_i\})+\sum_{i=k+1}^{m}\frac{h_{i}^{t}-h_i}{t}\mu_p((P^t)', \{u_i\})\!\!\bigg)\\
   &=& \frac{p-1}{n-p}\cdot \sum_{i=1}^{k}\mu_p(P'(z), \{u_i\}), 
\end{eqnarray*} 
{which is strictly positive by  Lemma \ref{relation between capacity measure and surface area measure}} and the fact that the origin $o$ is contained in at least one facet. It follows from the Minkowski inequality (\ref{Minkowski-classical-1}) that 
\begin{eqnarray*}
   0&<& \lim_{t\rightarrow 0^+}\frac{C_p((P^t)')-C_p((P^t)', P'(z))}{t}\\
   &\leq& \liminf_{t\rightarrow 0^+}\frac{C_p((P^t)')-C_p((P^t)')^{1-\frac{1}{n-p}}C_p(P'(z))^{\frac{1}{n-p}}}{t}\\
   &=& C_p(P'(z))^{1-\frac{1}{n-p}}\liminf_{t\rightarrow 0^+}\frac{C_p((P^t)')^{\frac{1}{n-p}}-C_p(P'(z))^{\frac{1}{n-p}}}{t}. 
\end{eqnarray*} 
Together with the fact that $P^t\subset P(z^t)$ for $t>0$ small enough, one has  $$C_p(P'(z))<C_p((P^t)')\leq C_p(P'(z^t))\ \ \ \mathrm{holds\ for\ some}\  t>0 \ \mathrm{small\ enough}.$$  This contradicts with the maximality of $C_p(P'(z))$ and hence $o\in P'(z)$.  Consequently  
$$
   C_p(P'(z))=\max_{x\in   M^*\cap\R^m_+}  C_p(P'(x)). 
$$ 
Lemma \ref{derivative of capacity} and the Lagrange multiplier rule yield  
\begin{eqnarray}\label{solution-Lagrange-1}
   (p-1)\cdot \mu_{p}(P'(z), \{u_i\}) &=&\eta\cdot \frac{\lambda_i}{\phi(z_i)} \ \ \ \ \mathrm{for\ \ all} \ \ \ 1\leq i\leq m, \\ \sum_{i=1}^{m}\lambda_i\varphi(z_i)&=&\varphi(1)\cdot \sum_{i=1}^{m} \lambda_i. \nonumber 
\end{eqnarray}  
Clearly $\eta>0$, as otherwise $\mu_{p}(P'(z), \{u_i\})=0$ and hence $S(P'(z), \{u_i\})=0$ (due to Lemma \ref{relation between capacity measure and surface area measure})  for all $1\leq i\leq m$. This leads to  the volume of $P'(z)$ equal to $0$,  which is impossible for a polytope $P(z)$ with nonempty interior. The positivity of $\eta$ further yields $\mu_p(P'(z), \{u_i\})>0$ and hence $z_i=h_{P(z)} (u_i)$ for all $1\leq i\leq m$. Moreover,  
$$
\sum_{i=1}^{m}\lambda_i\varphi(h_{P'(z)}(u_i))=\varphi(1)\cdot \sum_{i=1}^{m} \lambda_i \ \ \ \  \mathrm{and\ then} \ \ \ \  \|h_{P(z)}\|_{\varphi, \mu}=1.
$$

On the other hand, the constant $\eta$ satisfies the following formula 
\begin{eqnarray*} 
	\eta\cdot  \varphi(1)\cdot \sum_{i=1}^{m} \lambda_i&=&\eta \cdot \sum_{i=1}^{m}\big[\lambda_i \cdot \varphi(h_{P'(z)}(u_i))\big] \\&=&
   (p-1)\cdot \sum _{i=1}^m \big[\varphi(h_{P'(z)}(u_i))\cdot \phi(h_{P'(z)}(u_i))\cdot \mu_p(P'(z), \{u_i\})\big]\\ &=&(p-1)\cdot  \int_{\sphere} \varphi(h_{P'(z)}(u))\,d\mu_{p, \phi}(P'(z), u), 
\end{eqnarray*} 
where we have used  (\ref{solution-Lagrange-1}). 
Let   
\begin{eqnarray*} 
   \tau=\frac{p-1}{\eta}=\frac{\varphi(1) \cdot \int_{\sphere}\,d\mu(u)}{\int_{\sphere} \varphi(h_{P'(z)}(u))\,d\mu_{p, \phi}(P'(z), u)}, 
\end{eqnarray*}  
and then (\ref{solution-Lagrange-1}) yields $$\mu=\sum_{i=1}^{m}\lambda_i\delta_{u_i}=\tau\cdot \sum_{i=1}^m \Big[\phi(h_{P'(z)}(u_i)) \cdot \mu_p(P'(z), \{u_i\})\cdot \delta_{u_i}\Big].$$ Similarly, the constant $\tau$ can also be calculated by, due to (\ref{solution-Lagrange-1}), \begin{equation}\label{cal-tau-222} 
   \tau= \bigg(\frac{p-1}{n-p} \bigg)\cdot \frac{1}{C_p(P')} \cdot\int_{\sphere} \frac{h_{P'}(u)}{\phi(h_{P'}(u))}\,d\mu(u).  
\end{equation}  
This completes the proof.   \hfill $\Box$

\vskip 2mm When $\varphi=t^q/q$ with $q>1$, then $\phi=t^{1-q}$ and $(\phi, \varphi)$  satisfies conditions (A1) and (A2). In this case, the constant $\tau$ in Theorem \ref{discrete solution----1} can be calculated by  
\begin{eqnarray*} 
	\tau_q = \frac{ q\cdot \int_{\sphere}\,d\mu(u)}{\int_{\sphere}  h^{q}_{P'}(u)  \,d\mu_{p, q}(P', u)}=  \frac{c_q}{C_p(P')}, 
\end{eqnarray*}  
where the constant $c_q$ is
\begin{eqnarray}\label{constant c q}  
   c_q=\bigg(\frac{p-1}{n-p} \bigg)\cdot q  \cdot \int_{\sphere}   \,d\mu(u). 
\end{eqnarray}
Following immediately from  Theorem \ref{discrete solution----1}, one gets the solution for the normalized $p$-capacitary $L_q$ Minkowski problem for discrete measures.   The uniqueness is by Theorem \ref{uniqueness-mixed volume-1} or Corollary \ref{uniqueness-measure-1}. 
\bc \label{Lq--1--1-1-1-1}  
Let $\mu=\sum_{i=1}^{m}\lambda_i\delta_{u_i}$ be such that either $\mu(\{\xi\})=0$ or $\mu(\{-\xi\})=0$ for all $\xi\in S^{n-1}$, where $\lambda_1,\cdots, \lambda_m>0$ are given constants and $\{u_1, \cdots, u_m\} \subset S^{n-1}$ are not contained in any closed hemisphere. For $q>1$ and  $1<p<n$, the normalized $p$-capacitary $L_q$ Minkowski problem has a unique solution, i.e.,   there exists a unique polytope $P$ with the origin in its interior, such that,  
$$
   \mu=c_q \cdot \frac{\mu_{p, q}(P', \cdot)}{C_p(P')}.
$$ 
\ec

\subsection{The $p$-capacitary Orlicz-Minkowski problem for general measures}\label{M-capacity-general--111}

In this subsection, we provide a solution for the $p$-capacitary Orlicz-Minkowski problem for general measures. When $\varphi=t$, this has been investigated and solved in \cite{AGHLV, CNSXYZ, Jerison, Jerison-1996}.  See \cite[Theorem 1.2]{huang} for the volumetric case. We always use $\Omega$ to mean the interior of $\overline{\Omega}$. 

\bt\label{Orlicz-Minkowski-solution--1} 
Let $(\phi, \varphi)$ satisfy conditions (A1) and (A2).  Then the following are equivalent. 

\vskip 2mm \noindent i) $\mu$ is a nonzero finite Borel measure on $S^{n-1}$ whose support is not contained in any closed hemisphere, i.e.,  $$\int _{\sphere} \langle \eta, \theta\rangle_+ \,d\mu(\theta)>0\ \ \ \mathrm{for\ \  all\ \ } \eta\in S^{n-1};$$
ii)  There exist a constant $\tau>0$ and a convex body $\overline{\Omega}$ containing the origin $o$, such that, for $1<p<n$
$$
   \frac{\mu}{\phi(h_{\Omega})}=\tau \cdot \mu_p(\Omega,\cdot). 
$$ 
Moreover, if $\overline{\Omega}\in \cK_0$ is a convex body with the origin in its interior, then 
$$
   \mu =\tau\cdot \phi(h_{\Omega}) \cdot \mu_p(\Omega,\cdot)=\tau\cdot \mu_{p, \phi} (\Omega,\cdot). 
$$
\et 
\begin{proof}  
We first prove  $ i)\Rightarrow ii).$   Let $\mu$ be the given measure  satisfying with assumptions in Theorem \ref{Orlicz-Minkowski-solution--1}. Then  there exists a sequence of discrete measures $\mu_j$ defined on $S^{n-1}$ {satisfying that either $\mu_j(\{\xi\})=0$ or $\mu_j(\{-\xi\})=0$ for all $\xi\in S^{n-1}$ and} whose supports $\{u_1^j, \cdots, u_{m_j}^j\}$ are not contained in closed hemispheres, such that, $\mu_j \rightarrow \mu$ weakly as $j\rightarrow \infty$ (see, e.g., the proof of \cite[Theorem 7.1.2] {SchneiderBook}). By Theorem \ref{discrete solution----1},  there are polytopes $P_j$ with the origin in their interiors,  such that, for all $j\geq 1,$ 
\begin{equation}\label{multiple-discrete} 
    \frac{\mu_j}{ \phi(h_{P_j})} =\tau_j\cdot \mu_p(P'_j,\cdot)
\end{equation}  
with $\tau_j$ given by (\ref{cal-tau-222}) as follows: 
\begin{equation*} 
    \tau_j= \bigg(\frac{p-1}{n-p} \bigg)\cdot \frac{1}{C_p(P'_j)} \cdot\int_{\sphere} \frac{h_{P_j}(u)}{\phi(h_{P_j}(u))}\,d\mu_j(u).
\end{equation*} 
Moreover, by inequality (\ref{uniform-lower-bound}), $C_p(P'_j)\geq C_p(B^n_2)$ for all $j\geq 1$. 

 The radial function of a compact convex set $L\subset \Rn$, denoted by $\rho_L:\sphere\rightarrow [0, \infty)$, is defined by: for $u\in\sphere$,  $$\rho_L(u)=\max\{\lambda>0: \ \lambda u\in L\}.$$ For each $j=1, 2, \cdots,$  let $r_j=\max\{\rho_{P_j}(u): \ u\in \sphere\}$ be the maximal radius of $P_j$ and $v_j\in \sphere$ be a vector such that $r_j$ is obtained. Clearly the line segment $[0, r_jv_j]\subset P_j$ and hence  
$$
   r_j \cdot \langle u, v_j\rangle_+  \leq h_{P_j}(u) \ \ \ \mathrm{for\ all} \  \  \ u\in S^{n-1}.
$$ 
Note that $\|h_{P_j}\|_{\varphi, \mu_j}=1$ and  
$$
   \big\|r_j \cdot  \langle u, v_j\rangle_+ \big\|_{\varphi, \mu_j} =r_j\cdot\big\|  \langle u, v_j\rangle_+\big\|_{\varphi, \mu_j} \leq 1.
$$ 
A standard argument  (see e.g., \cite[Lemma 3]{HLYZ},  \cite[Corollary 3.7]{huang}, or similar results in \cite{ZHY2016}), as the supports of measures $\mu$ and $\mu_j$ are not contained in any closed hemisphere,  shows that there exists a constant $R>0$ such that $r_j\leq R$ for all $j\geq 1$, that is, $\{P_j\}_{j\geq 1}$ is bounded.  It follows from the Blaschke's selection theorem that there is a subsequence, which will not be relabeled, $\{P_j\}_{j\geq 1}$ converging to a compact convex set $\overline{\Omega}$ and $h_{P_j}\rightarrow h_{\overline{\Omega}}$ uniformly on $\sphere$. Moreover, $0<C_p(B^n_2)\leq C_p(\overline{\Omega})<\infty$ due to the continuity and monotonicity of the $p$-capacity and $\Omega\subset R\cdot B^n_2$. 

\vskip 2mm \noindent {\it Case 1: the interior of  $\overline{\Omega}$ is nonempty.} In this case, $\overline{\Omega}$  is a convex body containing the origin.  Let $\tau=\lim_{j\rightarrow \infty}\tau_j$.  By the continuity of the $p$-capacity,  \cite[Lemma 4.2]{ZHY2016},  and the uniform continuity of  the function $t/\phi(t)$  (whose value at $t=0$ is set to be $0$ due to $\lim_{t\rightarrow 0^+}t/\phi(t)=0$) on any closed bounded interval $[0, b]$, one has 
\begin{eqnarray*} 
    \tau  &=&\lim_{j\rightarrow \infty}  \bigg(\frac{p-1}{n-p} \bigg)\cdot \frac{1}{C_p(P'_j)} \cdot\int_{\sphere} \frac{h_{P_j}(u)}{\phi(h_{P_j}(u))}\,d\mu_j(u) \nonumber \\&=& \bigg(\frac{p-1}{n-p} \bigg)\cdot \frac{1}{C_p(\Omega)} \cdot\int_{\sphere} \frac{h_{\Omega}(u)}{\phi(h_{\Omega}(u))}\,d\mu(u).
\end{eqnarray*}  
Together with  (\ref{multiple-discrete}) and \cite[Lemma 4.1] {CNSXYZ},  one has $\mu_p(P'_j, \cdot)\rightarrow \mu_p(\Omega, \cdot)$ weakly and  for any continuous  function $f: \sphere\rightarrow \R$,   
\begin{eqnarray*}   
	\tau \cdot \int_{\sphere} f(u)\,d\mu_p(\Omega, u)&=& \lim_{j\rightarrow\infty} \tau_j\cdot \int_{\sphere} f(u)\,d\mu_p(P'_j, u)\\ &=& \lim_{j\rightarrow \infty} \int_{\sphere} \frac{f(u)}{ \phi(h_{P_j}(u))}\,d\mu_j(u)\\ &=&\int_{\sphere} \frac{f(u)}{ \phi(h_{\Omega}(u))}\,d\mu(u).
\end{eqnarray*} 
Hence, $\mu_p(P_j, \cdot)\rightarrow \frac{\mu}{\phi(h_{\Omega})}$ weakly and  $ 
\frac{\mu}{\phi(h_{\Omega})}=\tau \cdot \mu_p(\Omega,\cdot)  
$  by the uniqueness of the weak limit. Of course, if $\overline{\Omega}\in \cK_0$, then $h_{\Omega}$ is strictly positive on $\sphere$ and 
$$
  \mu =\tau\cdot \phi(h_{\Omega}) \cdot \mu_p(\Omega,\cdot)=\tau\cdot \mu_{p, \phi} (\Omega,\cdot). 
$$    
 
\noindent {\it Case 2: the interior of $\overline{\Omega}$ is empty.} In this case, without loss of generality, let 
$$
   \overline{\Omega}\subseteq\Big\{(x_1, \cdots, x_k, 0, \cdots, 0): \ \ x_1, \cdots, x_k\in \R\Big\}$$  
with $k$ the Hausdorff dimension of $\overline{\Omega}$ which is at most $n-1$. In fact $k>n-p$, as otherwise, $C_p(\overline{\Omega})=0$ which contradicts with $C_p(\overline{\Omega})>0.$

Recall that $h_{P_{j}}\rightarrow h_{\overline{\Omega}}$ uniformly on $\sphere$ and for all $j\geq 1$, $$h_{P_{j}}\le R\ \ \mathrm{on}\ \ \ \sphere,$$ where $R<\infty$ is the constant given above (i.e., the uniform upper bound of $r_j$). As $\phi$ is continuous and decreasing, one gets $\phi(h_{P_{j}}(u))\geq \phi(R)=: M>0$ for  all $u\in S^{n-1}$ and for all $j\geq 1$.  The constant $\tau_j$ can be bounded from below by (\ref{comp-phi-varphi-1}) as follows: 
\begin{eqnarray*}  \tau_j &=& \frac{\varphi(1) \cdot \int_{\sphere}\,d\mu_j(u)}{\int_{\sphere} \varphi(h_{P_j}(u))\,d\mu_{p, \phi}(P_j', u)}  \\ &\geq&  \frac{\varphi(1) \cdot \int_{\sphere}\,d\mu_j(u)}{\int_{\sphere}  h_{P_j}(u) \,d\mu_{p}(P_j', u)} \\ &=&  \frac{p-1}{n-p}\cdot \varphi(1) \cdot \frac{ \int_{\sphere}\,d\mu_j(u)}{C_p(P_j')}.
\end{eqnarray*} Hence, $\liminf_{j\rightarrow \infty} \tau_j \geq 2\tau_0$, if we let 
\begin{equation}\label{comp-phi-varphi-1}
2 \tau_0=\frac{p-1}{n-p}\cdot \varphi(1) \cdot \frac{ \int_{\sphere}\,d\mu(u)}{C_p(\overline{\Omega})}>0.
\end{equation}
Moreover,  (for convenience the Gauss maps of $P_{j}$ are all denoted by $\mathrm{g}$ unless otherwise stated)
\begin{eqnarray} \label{e5.46}
\int_{S^{n-1}}  \,d \mu&=&\lim_{j\to \infty} \int_{S^{n-1}}  \, d \mu_j \nonumber \\
&=& \lim_{j\to\infty}\tau_{j}\int_{\sphere} \phi(h_{P_{j}}(\theta)) \,d\mu_p(P'_{j},\theta)  \nonumber  \\
&\ge& \tau_0 M\cdot \liminf_{j\to \infty}\int_{\sphere}  \,d\mu_p(P'_{j},\theta).  
\end{eqnarray} On the other hand, one has \begin{equation}\label{lim is infty}\liminf_{j\to \infty}\int_{\sphere}  \,d\mu_p(P'_{j},\theta)=\infty,\end{equation} which was proved in \cite[Section 13.2]{AGHLV}. In fact,  (\ref{lim is infty}) follows directly from the combination of Propositions 13.5 and 13.6 in \cite{AGHLV} (by letting $f(y)=\|y\|^p$ for $y\in \Rn$) if $k=n-1$; while if $n-p<k<n-1$,  (\ref{lim is infty}) follows directly from \cite[(13.49)]{AGHLV} and its immediate consequence below (with $\mu_j$ in \cite[(13.49)]{AGHLV} replaced by $\mu_p(P_j', \cdot)$). Combining (\ref{e5.46}),  (\ref{lim is infty}), and the fact that $\mu$ is a finite measure on $S^{n-1}$, one gets a contradiction and hence the interior of $\overline{\Omega}$ cannot be empty. This completes the proof of $ i)\Rightarrow ii).$

Now we prove $ii)\Rightarrow i)$. Suppose that there exist a constant $\tau>0$ and a convex body $\overline{\Omega}$ containing the origin $o$, such that, for $1<p<n$
$$
   \frac{\mu}{\phi(h_{\Omega})}=\tau \cdot \mu_p(\Omega,\cdot). 
$$ 
Note that the support of $\mu_p(\Omega,\cdot)$ is  not contained in any closed hemisphere and 
$$
   \int_{\{\theta\in \sphere:\ h_{\Omega}(\theta)=0\}}\, d\mu_p(\Omega,\theta)=\int_{\{\theta\in \sphere:\ h_{\Omega}(\theta)=0\}}\frac{1}{\phi(h_{\Omega}(\theta))}\,d\mu(\theta)=0,
$$ 
where we have used $\frac{1}{\phi(t)}|_{t=0^+}=0$ due to $\phi(t)\rightarrow \infty$ as $t\rightarrow 0^+$.  Then,  for any given $\eta\in \sphere$, 
\begin{eqnarray*}
	0< \int_{\sphere} \langle \eta, \theta \rangle_+ \,d\mu_p(\Omega,\theta)=\lim_{k\rightarrow \infty} \int_{\{\theta\in \sphere:\ h_{\Omega}(\theta)\geq 1/k\}} \langle \eta, \theta \rangle_+ \,d\mu_p(\Omega,\theta).
\end{eqnarray*} 
Therefore, there exists $N_0$ (depending on $\eta$ of course) such that  for all $k\geq N_0$, 
\begin{eqnarray*}
	0<\int_{\{\theta\in \sphere:\ h_{\Omega}(\theta)\geq 1/k\}} \langle \eta, \theta \rangle_+ \,d\mu_p(\Omega, \theta).
\end{eqnarray*} 
This further implies that, for any $\eta\in \sphere$, 
\begin{eqnarray*}  
	\int_{\sphere} \langle \eta, \theta \rangle_+ \,d\mu(\theta)
	&\geq& \int_{\{\theta\in \sphere:\ h_{\Omega}(\theta)\geq  1/N_0\}} \langle \eta, \theta \rangle_+ \,d\mu(\theta)\\ 
	&=&\tau\cdot \int_{\{\theta\in \sphere:\ h_{\Omega}(\theta)\geq  1/N_0\}} \langle \eta, \theta \rangle_+\cdot \phi(h_{\Omega}(\theta)) \,d\mu_p(\Omega,\theta)\\ 
	&\geq&\tau\cdot \widetilde{m}_{\phi} \cdot \int_{\{\theta\in \sphere:\ h_{\Omega}(\theta)\geq  1/N_0\}} \langle \eta, \theta \rangle_+ \,d\mu_p(\Omega,\theta)\\ 
	&>& 0,
\end{eqnarray*} 
where we let $R=\sup_{\theta\in \sphere} h_{\Omega}(\theta)<\infty$ be a strictly positive constant and 
$$
   \widetilde{m}_{\phi}=\min\bigg\{\phi(t): \ \ \ t\in \Big[\frac{1}{N_0}, R\Big]\bigg\}\in (0, \infty). 
$$   
This completes the proof. 
\end{proof}

When $\varphi=t^q/q$ with $q>1$, then $\phi=t^{1-q}$ and $(\phi, \varphi)$  satisfies conditions (A1) and (A2). The solution for the normalized $p$-capacitary $L_q$ Minkowski problem for general measures can be stated as follows.  Let $c_q$ be the constant given by (\ref{constant c q}).

\bc  
Let $p\in (1, n)$ and $q>1$ be given constants. The following are equivalent. 

\vskip 2mm \noindent i)  $\mu$ is a nonzero finite Borel measure on $S^{n-1}$ whose support is not contained in any closed hemisphere,  i.e.,  
$$
   \int _{\sphere} \langle \eta, \theta\rangle_+ \,d\mu(\theta)>0\ \ \ \mathrm{for\ \  all\ \ } \eta\in S^{n-1};
$$  
ii) The normalized $p$-capacitary $L_q$ Minkowski problem has a unique solution, i.e.,   there exists a unique convex body $\overline{\Omega}$ containing the origin $o$, such that,  for $1<p<n$, 
$$
   \mu \cdot h^{q-1}_{\Omega}=c_q \cdot \frac{\mu_p(\Omega,\cdot)}{C_p(\Omega)}. 
$$ 
Moreover, if $\overline{\Omega} \in \cK_0$ is a convex body with the origin in its interior, then 
$$
   \mu  = \frac{c_q \cdot   \mu_{p, q} (\Omega,\cdot)}{C_p(\Omega)}. 
$$
\ec 
\begin{proof}  
The direction $ii)\Rightarrow i)$ follows immediately from Theorem \ref{Orlicz-Minkowski-solution--1} by letting $\phi=t^{1-q}$. For $i)\Rightarrow ii)$, the existence of a convex body $\overline{\Omega}$ containing the origin $o$ is an immediate consequence of Theorem \ref{Orlicz-Minkowski-solution--1}.  When $\overline{\Omega}\in \cK_0$, the uniqueness follows from Theorem \ref{uniqueness-mixed volume-1} or Corollary \ref{uniqueness-measure-1}. 
If the origin is not in the interior of $\overline{\Omega}$, the uniqueness can be proved based on the technique in proving  \cite[Lemma 2.1]{HLYZ}.  We include a self-contained brief proof here for completeness. 
  
Assume that there exist two convex bodies $\overline{\Omega}$ and $\overline{\Omega}_1$ containing the origin  such that 
$$
   \mu \cdot h^{q-1}_{\Omega} \cdot {C_p(\Omega)}=c_q \cdot  {\mu_p(\Omega,\cdot)} \ \ \ \mathrm{and}\ \ \  \mu \cdot h^{q-1}_{\Omega_1} \cdot {C_p(\Omega_1)}=c_q \cdot  {\mu_p(\Omega_1,\cdot)}. 
$$ 
These formulas yield that 
$$
   \bigg(\frac{p-1}{n-p}\bigg)\cdot  \int _{\sphere} h^q_{\Omega}\,d\mu(u)=\bigg(\frac{p-1}{n-p}\bigg)\cdot  \int _{\sphere} h^q_{\Omega_1}\,d\mu(u)=c_q.
$$
Let $\Sigma=\{u\in\sphere: h_{\Omega}(u)>0\}$ and hence  \begin{eqnarray*} 
	C_p(\Omega)&=&\bigg(\frac{p-1}{n-p}\bigg)\cdot \int_{\Sigma} h_{\Omega}\, d\mu_{p}(\Omega, u), \\  
	0 &=& C_p(\Omega) \int_{\sphere\setminus \Sigma}  h^{q-1} _{\Omega}(u) \,d\mu(u)= c_q \cdot \int_{\sphere\setminus \Sigma} \,d\mu_p(\Omega, u)\\ 
	C_p(\Omega, \Omega_1) &=&\bigg(\frac{p-1}{n-p}\bigg)\cdot \int_{\Sigma} h_{\Omega_1}\, d\mu_{p}(\Omega, u).
\end{eqnarray*}  
Holder's inequality implies that 
\begin{eqnarray*} 
	c_q =\bigg(\frac{p-1}{n-p}\bigg)\cdot  \int _{\sphere} h^q_{\Omega_1}\,d\mu(u) \geq \bigg(\frac{p-1}{n-p}\bigg)\cdot  \int _{\Sigma} h^q_{\Omega_1}\,d\mu(u) \geq c_q\cdot \bigg(\frac{C_p(\Omega, \Omega_1)}{C_p(\Omega)}\bigg)^q.  
\end{eqnarray*} 
Together with the Minkowski inequality (\ref{Minkowski-classical-1}), one gets $C_p(\Omega)\geq C_p(\Omega_1)$. By switching the roles of $\Omega$ and $\Omega_1$, one can also have $C_p(\Omega)\leq C_p(\Omega_1)$ and  $C_p(\Omega)= C_p(\Omega_1)$. Hence, the equality holds in Minkowski inequality (\ref{Minkowski-classical-1}) and then $\Omega_1$ is a translation of $\Omega$, say $\Omega_1=\Omega+a$ for some $a\in \Rn$. 

The uniqueness follows if $a=0$. To this end, assume that $a\neq 0$. By the translation invariance of the measures $\mu_p(\Omega, \cdot)$ and $\mu_p(\Omega_1, \cdot)$, one has 
$$
   \int_{\{u\in \sphere:\ \langle a, u\rangle >0\}} \big[ h^{q-1}_{\Omega_1}(u)-h^{q-1}_{\Omega}(u)\big]\,d\mu=0.
$$ 
Note that $h_{\Omega_1}(u)=h_{\Omega}(u)+\langle a, u\rangle>h_{\Omega}(u)$ for all $u\in \{u\in \sphere:\ \langle a, u\rangle >0\}$.  Then 
$$
   \int_{\{u\in \sphere:\ \langle a, u\rangle >0\}} \big[ h^{q-1}_{\Omega_1}(u)-h^{q-1}_{\Omega}(u)\big]\,d\mu>0,
$$ 
which follows from the assumption that the support of $\mu$ is not contained in the complement of $\{u\in \sphere: \langle a, u\rangle >0\}.$  This is a contradiction and hence $a=0$.  
\end{proof}

\vskip 2mm \noindent {\bf Acknowledgments.}  The second author is supported by a NSERC
grant. The authors are greatly indebted to Dr.\ J.\ Xiao for his valuable discussions and to the referee for many valuable comments which improve largely the quality of the paper.

\vskip 2mm \noindent Han Hong, \ \ \ {\small \tt honghan0917@126.com}\\
{ \em Department of Mathematics and Statistics,   Memorial University of Newfoundland,
   St.\ John's, Newfoundland, Canada A1C 5S7 }

\vskip 2mm \noindent Deping Ye, \ \ \ {\small \tt deping.ye@mun.ca}\\
{ \em Department of Mathematics and Statistics,
   Memorial University of Newfoundland,
   St.\ John's, Newfoundland, Canada A1C 5S7 }
   
   \vskip 2mm \noindent Ning Zhang, \ \ \ {\small \tt nzhang2@ualberta.ca}\\
{ \em Department of Mathematical and Statistical Sciences, University of Alberta, Edmonton, Alberta, Canada, T6G 2G1} 
   
\end{document}